\documentclass[final,oneside,onecolumn]{article}

\usepackage[body={6.4in,9in},top=0.9in,left=1in]{geometry}
\usepackage{overpic}
\usepackage{graphicx}
\usepackage{amsmath,amsthm,amssymb}
\usepackage{amsfonts}
\usepackage{paralist}
\usepackage{alltt}
\usepackage{url}
\usepackage{color}
\usepackage{algorithm,algorithmicx,algpseudocode}
\usepackage{amsmath}
\usepackage{amssymb}
\usepackage{enumitem}
\usepackage{mathtools}
\usepackage{cite}
\usepackage{caption}
\usepackage{subcaption}
\usepackage{colortbl}
\usepackage{tikz,cite,enumitem}
\usepackage{longtable}
\usepackage{longtable}
\usepackage{fancyhdr}
\usepackage{lineno}
\numberwithin{equation}{section}

\newcommand{\definedas}{\coloneqq}

\newcommand{\lon}{\phi}
\newcommand{\lat}{\theta}

\newcommand{\ds}{\displaystyle}

\usepackage{hyperref}

\numberwithin{figure}{section}

\numberwithin{table}{section}

\begin{document}

	
\title{Spectral interpolation in semi-implicit semi-Lagrangian methods for shallow water equations on the sphere}
	\author{Michael Chiwere\footnote{Malawi University of Business and Applied Sciences, Mathematical Sciences Department, Private Bag 303 Chichiri, Blantyre 3 (mchiwere@mubas.ac.mw).}  \and Daniel Fortunato\footnote{Center for Computational Mathematics \& Center for Computational Biology, Flatiron Institute (dfortunato@flatironinstitute.org).} \and Grady B. Wright \footnote{Boise State University, Department of Mathematics, 1910 University Drive, Boise, ID 83725-1555 (gradywright@boisestate.edu).}}
	\date{}
	\maketitle
	
\begin{abstract}
Semi-implicit semi-Lagrangian (SISL) methods are commonly used for the shallow water equations (SWE) because they allow for larger time steps than those permitted by the Courant-Friedrichs-Lewy (CFL) stability condition in Eulerian schemes. In these methods, the semi-Lagrangian treatment of advection is typically performed using lower-order interpolation, such as tensor-product Lagrange interpolation with cubic or quintic polynomials. However, operational SISL schemes routinely employ spectrally accurate spatial discretizations, such as spherical harmonics or the double Fourier sphere (DFS) method, for computing horizontal derivatives of the prognostic variables. This creates a mismatch in numerical accuracy, making the use of low-order interpolation less clearly justified.  In this work, we present the first numerical investigation of spectrally accurate interpolation in SISL schemes for the SWE. Our approach builds upon the recently developed DFS-based SWE model~\cite{yoshimura2022improved}, incorporating a spectral interpolation scheme that is accelerated using the nonuniform fast Fourier transform (NUFFT) to maintain the same overall computational complexity as the original model.  Using several standard SWE test cases, we evaluate the accuracy, conservation, and numerical diffusion of the new model, particularly over long integration times. Compared to an equivalent SISL model with low-order interpolation, the new model achieves higher accuracy, improved mass and energy conservation, and reduced numerical diffusion, demonstrating the potential benefits of incorporating spectrally accurate interpolation into SISL schemes.
\end{abstract}
	
	\section{Introduction\label{sec:intro}}
Explicit Eulerian schemes for global weather and climate models often require prohibitively small time steps to maintain stability due to the Courant-Friedrichs-Lewy (CFL) condition from the advection terms and the rapid propagation of gravitational waves.  
Semi-Lagrangian advection schemes partially alleviate this limitation, allowing larger time steps while preserving accuracy comparable to Eulerian methods~\cite{robert1981stable,staniforth1991semi}. However, the fast gravitational waves that arise from high-frequency fluctuations in the wind divergence still impose stability constraints on the size of the time step~\cite{robert1969integration,robert1972implicit}. Semi-implicit semi-Lagrangian (SISL) methods address this challenge by combining semi-Lagrangian advection with a semi-implicit treatment of terms that give rise to these gravity waves~\cite{robert1981stable,robert1982semi,yoshimura2012development}.  These methods are routinely used for the dynamical core of operational forecast models employed by weather forecasting centers throughout the world~\cite{tolstykh2017vorticity,ecmwf,figueroa2016brazilian,mengaldo2019current,yukimoto2019meteorological}.

A SISL scheme relies on four key components that influence its numerical properties~\cite{mengaldo2019current}: (i) the computation of departure point locations of fluid parcels in the semi-Lagrangian step, (ii) the temporal discretization of nonlinear terms responsible for gravity waves in the semi-implicit step, (iii) the spatial discretization of derivatives of the prognostic variables in the semi-implicit system, and (iv) the interpolation of the prognostic variables to the departure points. While all four components are important to the overall accuracy of the scheme, this study focuses on the latter two. High-order methods are commonly used for (iii), whereas relatively low-order methods are often employed for (iv).  However, beyond the motivation to reduce computational cost, there is limited numerical evidence in the literature to justify this discrepancy.  Notably, Giraldo et.\ al.~\cite{giraldo2003spectral} demonstrated the benefits of pairing high-order spectral element discretizations for both (iii) and (iv).  Building on this idea, we aim to further address this gap by investigating a new SISL scheme for the shallow water equations (SWE)---an idealized framework for studying the horizontal dynamics of weather and climate models---that employs global spectral spatial discretization for both (iii) and (iv).

\paragraph{Spectral spatial discretizations}
The spectral transform method, which represents solutions as spherical harmonic (SH) expansions, is commonly used for the horizontal spatial discretization in SISL schemes within many global operational weather/climate forecast models~\cite{ecmwf,figueroa2016brazilian,mengaldo2019current,yukimoto2019meteorological}. This approach, which we shorten to SHT (``spherical harmonic transform''), is widely used because it naturally resolves the singularity at the poles arising in spherical coordinates, provides exponential accuracy for approximating horizontal derivatives of (smooth) functions, and allows the Helmholtz equation that arises from the semi-implicit discretization to be solved efficiently in spectral space.  The main drawback of the SHT method, however, is the high computational cost associated with the Legendre transforms, which conventionally has a computational complexity of $\mathcal{O}(n^3)$, where $n$ is the truncated spectral wave number and the two-dimensional spherical grid employed has $\mathcal{O}(n^2)$ points.  Fast Legendre transforms have been developed that reduce the SHT cost to $\mathcal{O}(n^2 \log^2 n)$ (e.g.,~\cite{tygert2010fast,slevinsky2019fast}) and have been implemented in numerical and climate models~\cite{wedi2013fast}.  However, these methods can have high set-up costs and memory requirements.

The high computational cost of the SHT method originally led Merilees~\cite{merilees1973pseudospectral} to develop the Double Fourier Sphere (DFS) spectral method, which is based solely on Fourier series.  The key idea is to map functions defined on the sphere to a \textit{square} domain while preserving periodicity in both the azimuthal and polar directions.  This mapping results in an apparent `doubling' of the function in the polar direction that naturally facilitates an expansion in a bivariate Fourier series~\cite{townsend2016computing}, which explains the name DFS.  Compared to the SHT, DFS offers a promising alternative, as it relies solely on the fast Fourier transform (FFT), reducing the computational complexity to $\mathcal{O}\left(n^2\log n\right)$ with no setup costs and a small memory footprint.  Additionally, there are fast algorithms for solving the Poisson/Helmholtz equations with DFS~\cite{townsend2016computing,yee1981solution,cheong2000double,shen1999efficient}. 

Since Merilees' pioneering work, the DFS method has been employed in several numerical weather and climate models, as well as in a wide range of other applications. Layton \& Spotz were the first to develop a SISL scheme for the SWE using the DFS method~\cite{layton2003semi}.  However, their choice of Fourier basis produced discontinuities at the poles and non-isotropic waves, leading to potential instabilities. To mitigate this issue, they used an SHT filter of the prognostic variables, which increased the computational cost.  To eliminate the need for an SHT filter, Yoshimura and Matsumura adopted the Fourier basis introduced in~\cite{cheong2000double, cheong2000application} and developed a DFS-based SISL method that instead relies on higher-order diffusion for stability~\cite{yoshimura2005two,yoshimura2012development}, which could be implemented efficiently. More recently, Yoshimura introduced a new Fourier basis with a more stable method of computing the expansion coefficients and applied it to a SISL scheme for SWE~\cite{yoshimura2022improved}.  This approach eliminates the need for both an SHT filter and artificial diffusion while maintaining stability.

\paragraph{Common interpolation methods}
Although SHT and DFS spectral methods are widely used in SISL schemes, the interpolation step is typically handled using a completely different technique.  The most common choices are bivariate Lagrange interpolation with cubic or quintic polynomials, or bivariate cubic splines~\cite{cote1988two,ecmwf,layton2003semi,yoshimura2005two,yoshimura2012development,yoshimura2022improved}, which are relatively low order compared to the spectral accuracy of SHT and DFS methods. 
Two main arguments are typically given for using these lower-order interpolation methods. The first is rooted in the analysis of Falcone and Ferretti~\cite{FalconeFerrettiBook}, which establishes the error of semi-Lagrangian methods for the one-dimensional advection equation as $\mathcal{O}\left(\Delta t^k + \frac{h^p}{\Delta t}\right)$, where $k$ is the order of accuracy of the time-integration scheme, $p$ is the order of the spatial interpolation scheme, $h$ is the spatial resolution, and $\Delta t$ is the time step.  For most SISL methods, $k=2$, so choosing $p\geq 3$ ensures an error of $\mathcal{O}\left(\Delta t^2\right)$ when choosing $h=\mathcal{O}(\Delta t)$.  In practice for the SWE, cubic ($p=4$) or quintic ($p=6$) polynomial interpolation have been found to yield better results than the minimal quadratic interpolation ($p=3$) suggested by the theory for the simple advection equation~\cite{staniforth1991semi}. The second argument for using lower-order interpolation methods is computational efficiency. These methods have a complexity of $\mathcal{O}(n^2)$, whereas interpolation using SHT or DFS has a significantly higher complexity of $\mathcal{O}(n^4)$~\cite{layton2003semi}.

\paragraph{Motivation for spectral interpolation}
Despite the reasoning behind these arguments, we could not find any studies in the literature that numerically assess the benefits of 
employing spectral interpolation in SISL schemes.  This is an important question to explore, as spectral methods for advection-dominated problems are known to reduce diffusion and phase errors, better capture fine-scale structures and wave-like features, and better preserve accuracy over long integration periods~\cite{boyd2000chebyshev}.  The goal of this paper is to address this gap by developing the first SISL scheme for the SWE that employs spectral methods for both the spatial discretization of derivatives and interpolation to the departure points.  As part of this goal, we also aim to address the computational cost concern associated with spectral interpolation. 

The foundation for this scheme is the recent DFS-based SISL model developed by Yoshimura, for which an open-source implementation is available~\cite{yoshimura2022improved}.  This model introduces a new stable Fourier basis for DFS spatial discretization and employs cubic and quintic Lagrange interpolation. To ensure numerical consistency across all approximations, we replace its existing interpolation method with one also based on the DFS method using the same Fourier basis.  The main computational cost in this step is the evaluation of bivariate Fourier series of the prognostic variables at the departure points.  Na\"ively this step requires $\mathcal{O}(n^4)$ operations, but we significantly reduce this cost to $\mathcal{O}(n^2\log^2(1/\epsilon) + n^2\log n)$, where $\epsilon$ is the accuracy parameter, using a nonuniform Fast Fourier transform (NUFFT)~\cite{barnett2019parallel}.  This makes the computational complexity of the scheme consistent across all steps.  For further comparisons, we also develop an SHT-based SISL model, following~\cite{yukimoto2019meteorological}, that uses the SHT for the spatial discretization of derivatives and a similar DFS-based interpolation scheme.  

We evaluate these new schemes against existing SISL models that use lower-order interpolation methods on a standard suite of test cases for the SWE.  We examine the accuracy of the methods as the spatial and temporal steps are refined, their behaviors over long integration times, how well they conserve mass and energy, and the impact of numerical diffusion on the solutions.  Through these comparisons, we aim to determine whether spectral interpolation provides tangible advantages over conventional lower-order approaches.

\paragraph{Outline} The remainder of the paper is as follows.  In Section \ref{sec:ssl} we give an overview of the two-time-level SISL discretization. Next we briefly introduce the types of grids that are used in DFS and SHT methods in Section \ref{sec:grid}.  The DFS method using Yoshimura's stable Fourier basis is introduced in Section \ref{sec:newdfs} followed by a discussion about how the NUFFT can be used to reduce the computational complexity of DFS interpolation.  Numerical results for the new spectral SISL-DFS schemes are given in Section \ref{sec:res} together with a discusion.  Finally we make some concluding remarks in Section \ref{sec:conclusion}.

\section{Two-time-level semi-implicit semi-Lagrangian (SISL) scheme}
\label{sec:ssl}
SISL schemes are commonly implemented using either two or three-time levels. Proponents of two-time-level schemes argue that they are more efficient and easier to implement than their three-time-level counterparts and do not suffer from unstable temporal computational modes. However, as noted in~\cite{cote1988two}, the adoption of these schemes was initially hindered by difficulties in obtaining sufficiently accurate trajectories, which compromised the quality of the solution and rendered advantages of improved efficiency useless.

Temperton \& Staniforth~\cite{temperton1987efficient}, along with Donald \& Bates~\cite{mcdonald1987improving}, independently introduced the first-ever proper second-order methods for computing trajectories in two-time-level schemes. The former demonstrated that, for the shallow water finite-element model, the two-time-level approach was twice as efficient as the corresponding three-time-level scheme in~\cite{ritchie1988application}. The two-time-level scheme was further tested in a global spherical harmonic spectral model of the shallow water equations, showing a similar improvement in efficiency, being twice as fast as the three-time-level approach~\cite{cote1988two}.

In the context of DFS spectral shallow water models, three-time-level \cite{layton2003semi} and two-time-level \cite{yoshimura2022improved} schemes have been developed. Our study adopts the two-time-level scheme from~\cite{yoshimura2022improved} to discretize the rotating shallow water equations. However, our approach could successfully apply to three-time-level schemes and other two-time-level SISL schemes.

The SWE on a sphere can be expressed as follows~\cite{yoshimura2022improved} 
\begin{subequations}
\label{eq:swe}
	\begin{eqnarray}
		\frac{\mathrm{d}\left(\mathbf{v}+2\mathbf{\Omega}\times\mathbf{r}\right)}{\mathrm{d}t}	&=& -g\nabla h, \label{eq:swe1a}\\ 
		\frac{\mathrm{d}h}{\mathrm{d}t}	&=& -(h -h_s)\nabla \cdot \mathbf{v} + \mathbf{v}\cdot \nabla h_s, \label{eq:swe1b}
	\end{eqnarray}
\end{subequations}
where $t$ is time, $\mathbf{v}=(u,v)$ is the horizontal wind vector, $h$ represents the wave height, $h_s$ is the surface height or bathymetry, and $g$ is the gravitational acceleration, $\mathbf{\Omega}$ is the Earth's angular velocity vector, and $\mathbf{r}$ is the positional vector relative to the center of the Earth. These latter two terms combine to make up the Coriolis force $2\mathbf{\Omega}\times\mathbf{r}$. 

The SISL scheme used in~\cite{yoshimura2022improved} for \eqref{eq:swe} is a combination of the two-time-level SISL scheme~\cite{temperton1997treatment} and the Stable Extrapolation Two-Time-Level Scheme (SETTLS)~\cite{hortal2002development}.  Following the notation from~\cite{yoshimura2022improved}, the temporal discretization of the scheme  can be written as
\begin{subequations}
  \begin{equation}
	\begin{aligned}
		\frac{\left(\mathbf{v} +2\mathbf{\Omega}\times \mathbf{r}\right)^{+} - \left(\mathbf{v}+2\mathbf{\Omega}\times \mathbf{r}\right)_D^0}{\Delta t} &= -\frac{g\left(\nabla h_D^{(+)} + \nabla h^0\right)}{2} \\
		&+\beta_{\mathbf{v}}\frac{g\left(\nabla h_D^{(+)} + \nabla h^0\right)}{2} -\beta_{\mathbf{v}} \frac{g\left(\nabla h_D^0 + \nabla h^+\right)}{2},
    \end{aligned}\label{eq:swe2a}
  \end{equation}
  \begin{equation}
    \begin{aligned}
		\frac{h^+ - h_D^0}{\Delta t} &=-\frac{\left[\left(h - h_s\right)\nabla \cdot \mathbf{v}\right]_D^{(+)} + \left[(h-h_s)\nabla \cdot \mathbf{v}\right]^0}{2}  +\frac{\left[\mathbf{v}\cdot \nabla h_s\right]_D^{(+)} + \left[\mathbf{v}\cdot \nabla h_s\right]^0}{2} \\
		&+\beta_h\frac{\left[\bar{h}\nabla\cdot \mathbf{v}\right]_D^{(+)} + \left[\bar{h}\nabla \cdot \mathbf{v}\right]^0}{2} -\beta_h\frac{\left[\bar{h}\nabla \cdot\mathbf{v}\right]_D^0+\left[\bar{h}\nabla \cdot \mathbf{v}\right]^-}{2},
    \end{aligned}\label{eq:swe2b}
  \end{equation}
\end{subequations}
where $\Delta t$ is the time step. The notation is as follows:
\begin{itemize}
\item Terms in Roman font are scalars, while terms in boldface font are vectors.
\item Superscripts $-$, $0$, and $+$ denote the past ($t-\Delta t$), present ($t$), and future times ($t+\Delta t$), respectively, while the superscript $(+)$ refers to values extrapolated to the future time $t+\Delta t$. 
\item The subscript $D$ denotes the departure point, where the term is evaluated in space, while the absence of a subscript indicates evaluation of the term at the arrival point (see discussion below). 
\item The variable $\bar{h}$ is the homogeneous height used for linearizing \eqref{eq:swe1b}.
\item The parameters $\beta_v$ and $\beta_h$ are second-order ``decentering'' terms that affect the stability of the semi-implicit scheme~\cite{yukimoto2011meteorological}. 
Following~\cite{yoshimura2022improved}, we use the values $\beta_v=\beta_h=1$ in this study as they were found to be sufficient for ensuring stability in the DFS model.
\end{itemize}

A departure point $\mathbf{x}_D$ is defined as the upstream horizontal position along the trajectory in the $(\mathbf{x},t)$ plane of the fluid parcel that arrives at grid point $\mathbf{x}$ at a future time $t+\Delta t$ from the current time $t$~\cite{staniforth1991semi}. The calculation of the departure points is performed using an updated SETTLS scheme, ensuring stability for the time step sizes greater than 30 minutes as proposed by ~\cite{yoshimura2022improved}. The departure point is computed as follows:
\begin{subequations}
	\label{eq:settls}
	\begin{eqnarray}
		\mathbf{x}_D &=& \mathbf{x} - \frac{ \mathbf{v}_D^0 + \mathbf{v'}^+}{2}\Delta t, \qquad	\text{where, }	 \label{eq:nsettls} \\
		\mathbf{v'}^+ &=& \mathbf{v}_D^0 + \left(2\mathbf{\Omega }\times \mathbf{r}\right)_D -2\mathbf{\Omega}\times \mathbf{r} - \frac{g \left(\nabla h_D^{(+)} + \nabla h^0\right)}{2}\Delta t. \label{eq:esls}
	\end{eqnarray}
\end{subequations}

The provisional future velocity $\mathbf{v'}$ is obtained by solving Equation~\eqref {eq:swe1a} using an explicit semi-Lagrangian method. Simplifying \eqref{eq:settls} yields the following expression for the departure point:
\begin{equation}
	\label{eq:departurepoint}
	\mathbf{x}_D = \mathbf{x} - \Delta t \left[\left(\mathbf{v}^0 + \mathbf{\Omega}\times \mathbf{r} - \frac{g\Delta t \nabla h ^{(+)}}{4}\right)_D -
	\mathbf{\Omega} \times \mathbf{r} - \frac{g\Delta t \nabla h^0}{4}
	\right].
\end{equation}
Note that the value with subscript $D$ is computed at the departure point $\mathbf{x}_D$, but this term also depends on $\mathbf{x}_D$.  This nonlinear system is solved iteratively using fixed-point iteration~\cite{ritchie1996comparison,temperton2001two,layton2003semi}. 
Since the departure point $\mathbf{x}_D$ is generally not on the grid, the term with subscript $D$ in \eqref{eq:departurepoint} is computed via spatial interpolation from nearby grid points. 
We use quintic Lagrange interpolation instead of the cubic scheme employed by Yoshimura~\cite{yoshimura2022improved} for this sub-step, as we found that cubic interpolation was insufficient to preserve second-order accuracy in time.  Increasing the degree beyond five for this step did not lead to further improvements in accuracy\footnote{This differs from the interpolation required in \eqref{eq:sslforce}, where higher-order schemes do lead to lower errors and yield additional improvements, as demonstrated below.}.

Equations \eqref{eq:swe2a} and \eqref{eq:swe2b} can be simplified to
\begin{subequations}
    \label{eq:ssl}
    \begin{align}
		\mathbf{v}^+ + \frac{\beta_{\mathbf{v}} \Delta t}{2}g\nabla h^+ &= \mathbf{R}_\mathbf{v}, \label{eq:ssla} \\
        h^+ + \frac{\beta_h\Delta t}{2}\bar{h}\left(\nabla \cdot \mathbf{v}\right)^+ & = R_h, \label{eq:sslha} 
    \end{align}
\end{subequations}    
where
\begin{subequations}
\label{eq:sslforce}
\begin{align}
\mathbf{R}_\mathbf{v} &\definedas \mathcal{Q}\left[
\mathbf{v}^0 + 2\mathbf{\Omega}\times \mathbf{r} - \frac{\Delta t}{2}g\left((1-\beta_{\mathbf{v}}) \nabla h^{(+)}+\beta_{\mathbf{v}}\nabla h^0\right) 
\right]_D \!- 2\mathbf{\Omega} \times \mathbf{r} -\frac{\Delta t}{2}g(1-\beta_{\mathbf{v}})\nabla h^0 \label{eq:sslforcea}  \\
R_h &\definedas \left[h^0 + \frac{\Delta t}{2}\left[
		\left(\beta_h\bar{h}-\left(h-h_s\right)\right) (\nabla \cdot \mathbf{v})^{(+)}  + \mathbf{v} \cdot \nabla h_s^{(+)}
		-\beta_h \bar{h}\left(\nabla \cdot \mathbf{v}\right)^0 
		\right]\right]_D \nonumber \\
		&\,+\frac{\Delta t}{2}\left[
		\left(\beta_h\bar{h}-\left(h-h_s\right)\right)(\nabla \cdot \mathbf{v})^0 + \mathbf{v} \cdot  \nabla h_s^0\right].
\label{eq:sslforceb}
\end{align}
\end{subequations}
In the SISL scheme, the quantities $\mathbf{R}_\mathbf{v}$ and $R_h$ are transported along the trajectory to the arrival point ($A$). For the vector quantity $\mathbf{R}_\mathbf{v}$, it is necessary to account for the changing orientation of the coordinate system so that it is expressed in the local coordinate system, at the arrival point. This is achieved with the rotation matrix $\mathcal{Q}$ in \eqref{eq:sslforcea}, which is given explicitly by~\cite{temperton2001two}
\begin{equation}
\mathcal{Q} = 
    \begin{bmatrix}
        \phantom{-}p & q\\-q & p
    \end{bmatrix},
\end{equation}
where
\begin{align*}
   p &= \frac{\sin \lat_A\sin\lat_D + \left(1+\cos\lat_A\cos\lat_D\right)\cos (\lon_A - \lon_D)}{1 +\cos\eta}, \\
   q &= \frac{(\cos\lat_A + \cos\lat_D)\sin(\lon_A-\lon_D)}{1 +\cos\eta},
\end{align*}
with 
\[
p^2 + q^2 = 1.
\]
Here, $\phi$ and $\theta$ are longitude and colatitude, the subscripts $D$ and $A$ refer to the departure and arrival points respectively, and $\eta$ is the angle between radial position vectors $\mathbf{r}_D$ and $\mathbf{r}_A$, which is given as
\[
\eta = \cos ^{-1}\left( \frac{\mathbf{r}_A\cdot \mathbf{r}_D}
{|\mathbf{r}_A||\mathbf{r}_D|}
\right).
\]

Rather than solving for $\mathbf{v}^{+}$ directly in \eqref{eq:ssla}, the divergence and curl can be applied to this equation to obtain equations for the divergence $\delta^{+} = \nabla\cdot\mathbf{v}^{+}$ and vorticity $\zeta^{+} = \mathbf{k}\cdot\left(\nabla\times\mathbf{v}^{+}\right)$:
\begin{subequations}
\label{eq:divvort}
\begin{align}
\delta^{+} - \frac{\beta_{\mathbf{v}}\beta_h(\Delta t)^2}{4}g\overline{h} \nabla^2 \delta^{+} &= \nabla\cdot \mathbf{R}_{\mathbf{v}} - \frac{\beta_{\mathbf{v}}\Delta t}{2}g \nabla^2 R_h, \label{eq:divvorta}\\
\zeta^{+} &= \mathbf{k}\cdot\left(\nabla\times \mathbf{R}_{\mathbf{v}}\right), \label{eq:divvortb}
\end{align}
\end{subequations}
where $\mathbf{k}$ is the unit radial vector. Upon solving the Helmholtz equation \eqref{eq:divvorta} for the divergence and computing the vorticity \eqref{eq:divvortb}, the velocity $\mathbf{v}^+$ can be computed from
\begin{align}
    \mathbf{v}^+ = \nabla\chi^{+} -  \mathbf{k}\times \nabla{\psi^{+}},
\end{align}
where $\chi^{+}$ and $\psi^{+}$ are the velocity potential and the stream function, respectively, and are given as solutions to the following Poisson equations:
\begin{subequations}
\label{eq:vpsf}
\begin{align}
\nabla^2 \chi^{+} &= \delta^{+}, \\ 
\nabla^2 \psi^{+} &= \zeta^{+}.
\end{align}
\end{subequations}

The scheme developed in~\cite{yoshimura2022improved} uses the DFS method to expand the prognostic variables in \eqref{eq:ssl}--\eqref{eq:vpsf} in a Fourier basis with partial regularity on the sphere and computes the spatial derivatives in these equations in spectral space.  The Helmholtz \eqref{eq:divvorta} and Poisson \eqref{eq:vpsf} equations are also solved in Fourier space using a Galerkin method; see~\cite[Sec.\ 3.10]{yoshimura2022improved} for complete details.

While spectral discretizations are used to compute the derivatives in this SISL model, quantities that are evaluated at the off-grid departure points (indicated by the subscript $D$) in \eqref{eq:sslforce} are computed using Lagrange polynomial interpolation (specifically, using quintic polynomials for the quantity in \eqref{eq:sslforcea} and cubic polynomials for \eqref{eq:sslforceb}). Note that this is common practice in many other SISL spectral models~\cite{cote1988two,ecmwf,layton2003semi,yoshimura2005two,yoshimura2012development,yukimoto2019meteorological}. In this study, we propose replacing the Lagrange interpolation technique with one based on the DFS method (see \cite{chiwere2025barycentric}). The technique uses the same Fourier basis as used for the spatial discretization in the DFS SISL model of~\cite{yoshimura2022improved}.  We also adapt this interpolation method for the SHT SISL model from~\cite{yukimoto2019meteorological}.

    \section{Grid types\label{sec:grid}}
In this study, we consider four different tensor-product latitude-longitude grids for the spatial discretizations. Each one uses equally spaced points in the longitudinal (zonal) direction, but different arrangements of grid points in the latitudinal (meridional) direction. The first three are used in the DFS SISL model and are described in~\cite{yoshimura2022improved} as follows:
\begin{subequations}
\label{eq:dfs_swe_grids}
	\begin{alignat}{4}
			\text{Grid [-1]:}\quad \theta_j &= \frac{j}{J}\pi, &&\quad j=1,\dots,J-1, &&\quad N_{\rm lat}=J-1,\\
			\text{Grid [0]:}\quad \theta_j &= \frac{j+1/2}{J}\pi, &&\quad j=0,\dots,J-1, &&\quad N_{\rm lat}=J,\\
			\text{Grid [1]:}\quad \theta_j &= \frac{j}{J}\pi, &&\quad j=0,\dots,J, &&\quad N_{\rm lat}=J+1.
		\end{alignat}
\end{subequations}
Here, $\theta_j$ denotes the colatitude at each grid point, and $J$ is an integer related to the number of latitudinal grid points $N_{\rm lat}$ chosen.  All grids use $M_{\rm lon}=2J$ equally spaced points in longitude. Grid [0] is the most commonly used grid type in geophysical problems~\cite{merilees1973pseudospectral,cheong2000double,cheong2000application,yoshimura2012development,yoshimura2022improved}. For all these grids, the DFS method from~\cite{yoshimura2022improved} ensures smooth approximations of functions and their derivatives at the poles (see Section \ref{sec:newdfs}).

The fourth grid type is used in the SHT SISL model and employs Gauss--Legendre (GL) points in latitude:
\begin{align}
    \text{GL:}\quad\theta_j&=\arccos({z_j}), \quad j=0, \dots , J-1, \quad N_{\rm lat}=J, \label{eq:GLGrid}
\end{align}
where $ z_j$ are roots of the degree-$J$ Legendre polynomial. Note that this is the same GL grid in~\cite{chiwere2025barycentric}, where we used it to show that it is possible to develop a DFS interpolation method.  The GL grid is most commonly used for problems involving the SHT. 

\section{Partial regularity DFS expansion method\label{sec:newdfs}}
We review the DFS approach for computing the discrete bivariate Fourier series on the sphere from~\cite{yoshimura2022improved} that forms the foundation for our interpolation method. This approach uses a new basis for the expansion and least squares projection for computing the corresponding expansion coefficients. The new basis enforces partial regularity, ensuring continuity of both scalar and vector fields and their derivatives at the poles, which is not guaranteed by alternative bivariate Fourier methods. Continuity of the vector field at the poles facilitates direct interpolation of the wind velocity components $u$ and $v$ without the need to adjust them to $u\sin\theta$ and $v\sin\theta$, as other approaches require.

The method approximates a scalar real-valued function $f(\phi, \theta)$ on the sphere using the expansion
\begin{align}
    f^{N,M}(\phi,\theta) = \sum_{m=0}^M f_{m}^{{\rm c}, N}(\theta)\cos(m\lon) + \sum_{m=1}^M f_{m}^{{\rm s}, N}(\theta)\sin(m\lon),
    \label{eq:pr_dfs}
\end{align}
where
\begin{align}
    f_{m}^{{\rm c(s)}, N} = 
\begin{cases}
    \ds \sum_{n=0}^N f_{n,0}^{{\rm c},N} \cos(n\theta), & m = 0, \\
    \ds \sum_{n=0}^{N-1} f_{n,1}^{{\rm c(s)},N} \sin(\theta)\cos(n\theta), & m = 1, \\
    \ds \sum_{n=1}^{N-1} f_{n,m}^{{\rm c(s)},N} \sin(\theta)\sin(n\theta), & m\geq 2\;\text{and even,}\\   
    \ds \sum_{n=1}^{N-2} f_{n,m}^{{\rm c(s)},N} \sin^2(\theta)\sin(n\theta), & m\geq 3\;\text{and odd.}
\end{cases}
\end{align}
Here, $\phi$ and $\theta$ are the azimuthal and polar coordinates, respectively, and the superscripts c and s on the coefficients correspond to the azimuthal cosine and sine expansions, respectively. Before discussing how the expansion coefficients are computed, we note that: (1) for each even $m$, $f^{N,M}$ is $\pi$-periodic in $\lon$ and even in $\theta$; (2) for each odd $m$, $f^{N,M}$ is $\pi$-anti-periodic in $\lon$ and odd in $\theta$; and (3) $f^{N,M}$ is constant when $\theta=0,\pi$.  Thus, by Lemma 1 in~\cite{chiwere2025barycentric} it is a continuous function on the sphere.  Furthermore, it can be shown that both components of the gradient in spherical coordinates, i.e., $\nabla f^{N,M} = \frac{1}{\sin\theta}\frac{df^{N,m}}{d\phi}\boldsymbol{\lon} + \frac{df^{N,m}}{d\theta}\boldsymbol{\theta}$, are well-behaved at the poles~\cite[Section 3.2]{yoshimura2022improved}.  These properties are why we call \eqref{eq:pr_dfs} the partial regularity DFS expansion; full regularity would ensure smoothness at the poles for all derivatives, as is the case for a spherical harmonic expansion.


The expansion coefficients in \eqref{eq:pr_dfs} are determined using a least squares procedure that involves first expanding $f$ in the standard DFS basis, which we express as
\begin{align}
\overline{f}^{N,M}(\phi,\theta) = \sum_{m=0}^M \overline{f}_{m}^{{\rm c}, N}(\theta)\cos(m\lon) + \sum_{m=1}^M \overline{f}_{m}^{{\rm s}, N}(\theta)\sin(m\lon),
\end{align}
where
\begin{align}
\overline{f}^{{\rm c(s)},N} = 
\begin{cases}
\ds \sum_{n=0}^N \overline{f}_{n,m}^{{\rm c(s)}}\cos(n\theta), & \text{$m$ even,} \\
\ds \sum_{n=1}^N \overline{f}_{n,m}^{{\rm c(s)}}\sin(n\theta), & \text{$m$ odd.}
\end{cases}
\end{align}
The coefficients $\overline{f}_{n,m}^{{\rm c(s)}}$ are computed from samples of $f$ on the grids \eqref{eq:dfs_swe_grids} using the FFT of $f$ in $\phi$ and the fast cosine and sine transforms (FCT and FST, respectively) in $\theta$. The expansion coefficients in \eqref{eq:pr_dfs} are then determined by picking them to minimize the square of the $L_2$-norm of the residual $f^{N,M}(\phi,\theta) - \overline{f}^{N,M}(\phi,\theta)$, i.e., minimizing
\begin{align}
    E_2(f_{n,m}^{{\rm c(s)}}) = \frac{1}{2\pi^2}\int_{0}^{2\pi}\int_{0}^{\pi} \left[f^{N,M}(\phi,\theta) - \overline{f}^{N,M}(\phi,\theta)\right]^2 d\theta d\phi.
\end{align}
This can be viewed as an $L_2([0,2\pi]\times[0,\pi])$-projection of the expansion of $f$ in the standard DFS basis onto the partial regularity DFS basis.
As derived in~\cite{yoshimura2022improved}, this leads to the following conditions:
\begin{align*}
    \int_{0}^{\pi} \left(f^{{\rm c},N}_0(\theta) - \overline{f}^{{\rm c},N}_0(\theta)\right)\cos(n\theta)d\theta &= 0, \\
    \int_{0}^{\pi} \left(f^{{\rm c(s)},N}_1(\theta) - \overline{f}^{{\rm c(s)},N}_1(\theta)\right)\sin(\theta)\cos(n\theta)d\theta &= 0, \\
    \int_{0}^{\pi} \left(f^{{\rm c(s)},N}_m(\theta) - \overline{f}^{{\rm c(s)},N}_m(\theta)\right)\sin(\theta)\sin(n\theta)d\theta &= 0, \; \text{for $m\geq 2$ and even,} \\    
    \int_{0}^{\pi} \left(f^{{\rm c(s)},N}_m(\theta)(\theta) - \overline{f}^{{\rm c(s)},N}_m(\theta)\right)\sin^2(\theta)\sin(n\theta)d\theta &= 0, \; \text{for $m\geq 3$ and odd.}
\end{align*}
For the different values of $m$ in the above equations, linear systems of equations can be worked out relating the coefficients $\overline{f}^{{\rm c(s)},N}_{n,m}$ to $f^{{\rm c(s)},N}_{n,m}$. 
For $m=0$ there is a direct relationship, while for $m=1$ and even $m\geq 2$ one gets tridiagonal systems to solve and for odd $m\geq 3$ one gets pentadiagonal systems; see~\cite[Section 3.3]{yoshimura2022improved} for full details on these systems.  The cost of solving these systems for all $m$ is $\mathcal{O}(MN)$.

Once the coefficients $f_m^{{\rm c(s)},N}$ have been computed, an inverse fast transform is needed to compute $f^{N,M}$ in \eqref{eq:pr_dfs} at the grid points.  This can be accomplished by converting $f^{N,M}$ to the standard DFS basis 
\begin{align}
\label{eq:std_dfs}
f^{N,M}(\phi,\theta) = \sum_{m=0}^M \tilde{f}_{m}^{{\rm c}, N}(\theta)\cos(m\lon) + \sum_{m=1}^M \tilde{f}_{m}^{{\rm s}, N}(\theta)\sin(m\lon),
\end{align}
with 
\begin{align}
\tilde{f}^{{\rm c(s)},N} = 
\begin{cases}
\ds \sum_{n=0}^N \tilde{f}_{n,m}^{{\rm c(s)}}\cos(n\theta), & \text{$m$ even,} \\
\ds \sum_{n=1}^N \tilde{f}_{n,m}^{{\rm c(s)}}\sin(n\theta), & \text{$m$ odd,}
\end{cases}
\end{align}
for which the inverse FFT, FCT, and FST can be used. 
The relationship between the coefficients in \eqref{eq:pr_dfs} and \eqref{eq:std_dfs} is given as follows~\cite[Section 3.1]{yoshimura2022improved}:
\begin{itemize}
\item For $m=0$:
\begin{align*}
\tilde{f}_{n,0}^{\rm c} &= f_{n,0}^{\rm c},\; \text{for $n=0,\ldots,N$.}
\end{align*}

\item For $m=1$:
\begin{align*}
\tilde{f}_{1,1}^{\rm c(s)} &= \frac{1}{2}\left(2f_{0,1}^{\rm c(s)} - f_{2,1}^{\rm c(s)}\right),\; \\
\tilde{f}_{n,1}^{\rm c(s)} &= \frac{1}{2}\left(f_{n-1,1}^{\rm c(s)} - f_{n+1,1}^{\rm c(s)}\right),\;  \text{$n=2,\ldots,N-2$,} \\
\tilde{f}_{n,1}^{\rm c(s)} &= \frac{1}{2}f_{n-1,1}^{\rm c(s)},\; \text{$n=N-1,N$}. 
\end{align*}

\item For even $m\geq 2$:
\begin{align*}
\tilde{f}_{0,m}^{\rm c(s)} &= 0,\; \tilde{f}_{1,m}^{\rm c(s)} = \frac{1}{2}f_{2,m}^{\rm c(s)},\; \\
\tilde{f}_{n,m}^{\rm c(s)} &= \frac{1}{2}\left(-f_{n-1,m}^{\rm c(s)} + f_{n+1,m}^{\rm c(s)}\right),\;  \text{$n=2,\ldots,N-2$,} \\
\tilde{f}_{n,m}^{\rm c(s)} &= -\frac{1}{2}f_{n-1,m}^{\rm c(s)},\; \text{$n=N-1,N$}. 
\end{align*}

\item For odd $m\geq 3$:
\begin{align*}
\tilde{f}_{0,m}^{\rm c(s)} &= 0,\; \tilde{f}_{1,m}^{\rm c(s)} = \frac{1}{4}\left(3f_{1,m}^{\rm c(s)} - f_{3,m}^{\rm c(s)}\right),\; \tilde{f}_{2,m}^{\rm c(s)} = \frac{1}{4}\left(2f_{2,m}^{\rm c(s)} - f_{4,m}^{\rm c(s)}\right),\; \\
\tilde{f}_{n,m}^{\rm c(s)} &= \frac{1}{4}\left(-f_{n-2,m}^{\rm c(s)} + 2f_{n,m}^{\rm c(s)} - f_{n+2,m}^{\rm c(s)}\right),\;  \text{$n=3,\ldots,N-3$,} \\
\tilde{f}_{n,m}^{\rm c(s)} &= \frac{1}{4}\left(-f_{n-2,m}^{\rm c(s)} + 2f_{n,m}^{\rm c(s)}\right),\; \text{$n=N-2,N-1$},\;  \tilde{f}_{N,m}^{\rm c(s)} = -\frac{1}{4}f_{N-2,m}^{\rm c(s)}.
\end{align*}
\end{itemize}


To expand vector fields $\mathbf{v} = (u,v)$ on the sphere,~\cite{yoshimura2022improved} starts by considering the Helmholtz decomposition of the field: 
\begin{align}
   \label{eq:windpotentialstream}
   u &= \dfrac{1}{a\sin\theta}\dfrac{\partial \chi}{\partial \phi} +\dfrac{1}{a}\dfrac{\partial \psi}{\partial \theta}, \\
   v &= -\dfrac{1}{a}\dfrac{\partial \chi}{\partial \theta} + \dfrac{1}{a\sin\theta}\dfrac{\partial \psi}{\partial \phi},
\end{align}
where $a$ is the Earth's radius, and the scalar potentials $\chi$ and $\psi$ are the velocity potential and stream function. The scalar potentials are then expanded into the partial regularity DFS basis similar to~\eqref{eq:pr_dfs}.  A least squares procedure is again used to compute the coefficients of these expansions from the standard DFS series of the components of $u$ and $v$.  Once these coefficients are determined, the partial regularity DFS series for each component are expressed in the standard DFS series for which the FFT, FCT, and FST can be used.  We omit these details and refer the reader to~\cite[Section 3.5--3.6]{yoshimura2022improved}.

\subsection{Zonal filter}	
\label{sec:zonalfilter}
In a typical tensor-product latitude-longitude grid, the spacing between zonal grid points becomes progressively narrower as we approach the poles. This narrowing can lead to computational instabilities, especially when using larger time steps in polar regions. A standard solution to this issue in DFS models is the zonal filter~\cite{merilees1974numerical,boer1975fourier,cheong2000double,yoshimura2022improved}, which filters out high zonal wavenumber components in the polar regions, reducing oversampling at the poles and ensuring a more uniform resolution over the sphere. An alternative idea to achieve the same goal is to use a reduced grid in value space (e.g.~\cite{hortal1991use,malardel2016new}). 

In this work, we use the same zonal filter from~\cite{yoshimura2022improved} and set the truncation wavenumber in longitude in the expansion \eqref{eq:pr_dfs} as:
\begin{equation}
    M_{\rm zf}(\theta_j) = \min (M, M_0+M\sin(\theta_j)),
\end{equation}
where $M_0=20$. In this formulation, all expansion coefficients $f_{n,m}^{\rm c(s)}$ for the wavenumber $m>M_{\rm zf}(\theta )$ are set to zero, effectively filtering out high wavenumber components at each latitude.

\section{Fast DFS interpolation method\label{sec:fast_dfs_interp}}
The idea of our new proposed spectral interpolation method for the DFS SISL scheme is to expand the terms in $\mathbf{R}_{\mathbf{v}}$ and $R_h$ in \eqref{eq:sslforce} with the subscript $D$, namely, 
\begin{align}
\mathbf{P}_\mathbf{v} &= 
\mathbf{v}^0 + 2\mathbf{\Omega}\times \mathbf{r} - \frac{\Delta t}{2}g\left((1-\beta_{\mathbf{v}}) \nabla h^{(+)}+\beta_{\mathbf{v}}\nabla h^0\right), \label{eq:deptforcea}  \\
P_h &= h^0 + \frac{\Delta t}{2}\left[
		\left(\beta_h\bar{h}-\left(h-h_s\right)\right) (\nabla \cdot \mathbf{v})^{(+)}  + \mathbf{v} \cdot \nabla h_s^{(+)}
		-\beta_h \bar{h}\left(\nabla \cdot \mathbf{v}\right)^0 
		\right],
\label{eq:deptforceb}
\end{align}
using the vector and scalar partial regularity DFS techniques discussed in the previous section and evaluate them at the departure points $\mathbf{x}_{D}$. For an $N_{\rm lat}\times M_{\rm lon}$ grid, the number of departure points is $M_{\rm lon} N_{\rm lat}$, so the cost of evaluating the DFS expansions directly is $\mathcal{O}((M_{\rm lon} N_{\rm lat})^2)$, assuming the series truncation parameters satisfy $M=\mathcal{O}(M_{\rm lon})$ and $N=\mathcal{O}(N_{\rm lat})$. This is significantly more expensive than the $\mathcal{O}(M_{\rm lon} N_{\rm lat} \log(M_{\rm lon} N_{\rm lat}))$ cost associated with the rest of the DFS SISL scheme.  

However, these series can be recast as bivariate Fourier series in complex exponential form.  Specifically, if the DFS series for $P_h$ or the components of $\mathbf{P}_{\mathbf{v}}$ have expansions in the form of \eqref{eq:std_dfs}, they can be equivalently written as
\begin{align}
\label{eq:std_dfs_exp}
f^{N,M}(\phi,\theta) = \sum_{m=-M}^{M} 
\sum_{n=-N}^{N}\tilde{f}_{n,m} e^{i m \phi}e^{i n\theta},
\end{align}
where 
\begin{gather}
\label{eq:dfs_fourier_coeffs}
\begin{aligned}
\tilde{f}_{0,0} &= \tilde{f}_{0,0}^{\rm c},\; \tilde{f}_{0,m} = 0\; \text{for odd $m\geq 1$},\; \tilde{f}_{0,m} = \frac12 \tilde{f}_{0,m}^{\rm c}\; \text{for even $m\geq 2$},\\
\tilde{f}_{\pm n,m} &= \frac14 \left(\tilde{f}_{n,m}^{\rm c} - i\tilde{f}_{n,m}^{\rm s}\right),\; 
\tilde{f}_{\pm n,-m} = \frac14 \left(\tilde{f}_{n,m}^{\rm c} + i\tilde{f}_{n,m}^{\rm s}\right)\; \text{for even $m\geq 2$ and $n\geq 1$},\\
\tilde{f}_{\pm n,m} &= \mp \frac14 \left(i\tilde{f}_{n,m}^{\rm c} + \tilde{f}_{n,m}^{\rm s}\right),\;
\tilde{f}_{\pm n,-m} = \mp \frac14 \left(i\tilde{f}_{n,m}^{\rm c} - \tilde{f}_{n,m}^{\rm s}\right)\; \text{for odd $m\geq 1$ and $n\geq 1$}.
\end{aligned}
\end{gather}
The relationships between the coefficients follow from~\cite[Section 2]{wright2015extension}.  Evaluating \eqref{eq:std_dfs_exp} at $M_{\rm lon} N_{\rm lat}$ points corresponds to a type-II nonuniform discrete Fourier transform (NUDFT)~\cite{greengard2004accelerating}.  

Several NUFFT algorithms have been developed for computing  type-II NUDFTs (as well as other types) in near optimal complexity (e.g.~\cite{dutt1993fast,greengard2004accelerating,potts2003fast,ruiz2018nonuniform,barnett2019parallel}), and there are many open-source implementations available (e.g., \cite{lin2018python,ruiz2018nonuniform,KeKuPo09,Slevinsky2025}). In this work we employ the NUFFT from~\cite{barnett2019parallel}, which evaluates \eqref{eq:std_dfs_exp} in  $\mathcal{O}\left(M_{\rm lon} N_{\rm lat}\log^2(1/\epsilon) + M N\log(M N)\right)$ operations, where $\epsilon$ controls the desired accuracy.  We use the implementation of this algorithm in the FINUFFT library~\cite{Barnett2025}, which supports fast algorithms for NUDFTs of types I--III in one, two, and three dimensions. This library also provides a Fortran interface that is particularly useful here, since the SWE code~\cite{yoshimura2022improved} utilized in this study is written in Fortran.

\subsection{Interpolation for the SHT SISL scheme}	
We adopt a similar strategy for interpolation in the SHT SISL model from~\cite{yukimoto2019meteorological}.  For this model, the SHT is used to obtain the values of \eqref{eq:deptforcea} and \eqref{eq:deptforceb} at an $M_{\rm lon}\times N_{\rm lat}$ tensor-product grid, with the latitude points arranged on a GL grid \eqref{eq:GLGrid}.  Following~\cite{chiwere2025barycentric}, the DFS method can be used to produce bivariate Fourier interpolants from these values in the form of~\eqref{eq:std_dfs_exp}.  The Fourier coefficients of ~\eqref{eq:std_dfs_exp} can be computed in the longitude direction using the FFT, but the nonequispaced GL grid restricts its use in the latitude direction.  For this computation, we use a direct inversion of the DFT matrix, which leads to a computational cost for computing the Fourier coefficients of  
$\mathcal{O}(M_{\rm lon}N_{\rm lat}^2 + M_{\rm lon}N_{\rm lat}\log M_{\rm lon})$, after pre-computing an LU factorization of the DFT matrix.  Once the Fourier coefficients for \eqref{eq:deptforcea} and \eqref{eq:deptforceb} are computed,  we use the NUFFT to evaluate the interpolants at the departure points, similar to the technique described above.

We make a few comments about this interpolation approach:
\begin{itemize}
     \item The cost of computing the Fourier coefficients in latitude can be improved using fast inverse NUFFT methods (e.g.~\cite{ruiz2018nonuniform,wilber2025}), but we did not employ these here.
     \item We did not find it necessary to use the partial regularity DFS basis for expanding \eqref{eq:deptforcea} and \eqref{eq:deptforceb} to maintain temporal stability during the simulations. This may be because this SISL model continually employs the forward and inverse SHT of the prognostic variables, which has been noted to produce stabilizing effects in other DFS-based SISL methods~\cite{layton2003semi}.
     \item The spectrally accurate SHT interpolation technique from~\cite{belkner2024cunusht}, which is also based on the DFS method, could also be used as an alternative to the one we employ to achieve better speed-ups.
\end{itemize}

\section{Numerical results \label{sec:res}}
 We evaluate the DFS SISL scheme with the new spectral interpolation method using several standard shallow water test cases, comparing it against the low order Lagrange interpolation approach in Yoshimura~\cite{yoshimura2022improved}. The test cases include Williamson test cases 1, 2, 5, and 6~\cite{williamson1992standard} and the Galewsky test case~\cite{galewsky2004initial}. For reference, we denote the results obtained using~\cite{yoshimura2022improved} as \textit{LAG Interp} and the ones obtained from the new proposed method as \textit{DFS Interp}.  We also include results for the SHT SISL scheme on some of these test cases, and these are also denoted by \textit{DFS Interp} with a GL grid.

The results of the DFS SISL schemes are presented using the $N_{\rm lat}\times M_{\rm lon}$ common grids [-1], [0], and [1] discussed in Section \ref{sec:grid}.  These grids are determined by the parameter $J$, where $M_{\rm lon} = 2J$ and $N_{\rm lat} = J+[\text{grid number}]$.  The SHT SISL scheme uses a GL grid with $M_{\rm lon} = 2J$ and $N_{\rm lat} = J$.

	
\subsection{Williamson Test Case 1 (TC1)\label{subsec:tc1}}
This test case considers the advection of a compactly supported cosine bell using a constant wind field that produces solid body rotation at an angle $\alpha=\pi/2-0.05$ with respect to the polar axis. The rotation period is 12 days, so the unaltered bell returns to its original position at this time, where errors can be measured. Since the test case is only for pure advection of a scalar height field, the momentum equations do not need to be solved, and the test case really only involves computing departure points and interpolation of the height field.  

\begin{figure}[tb]
\centering
\begin{subfigure}[b]{0.45\textwidth}
\centering
\includegraphics[width=\textwidth]{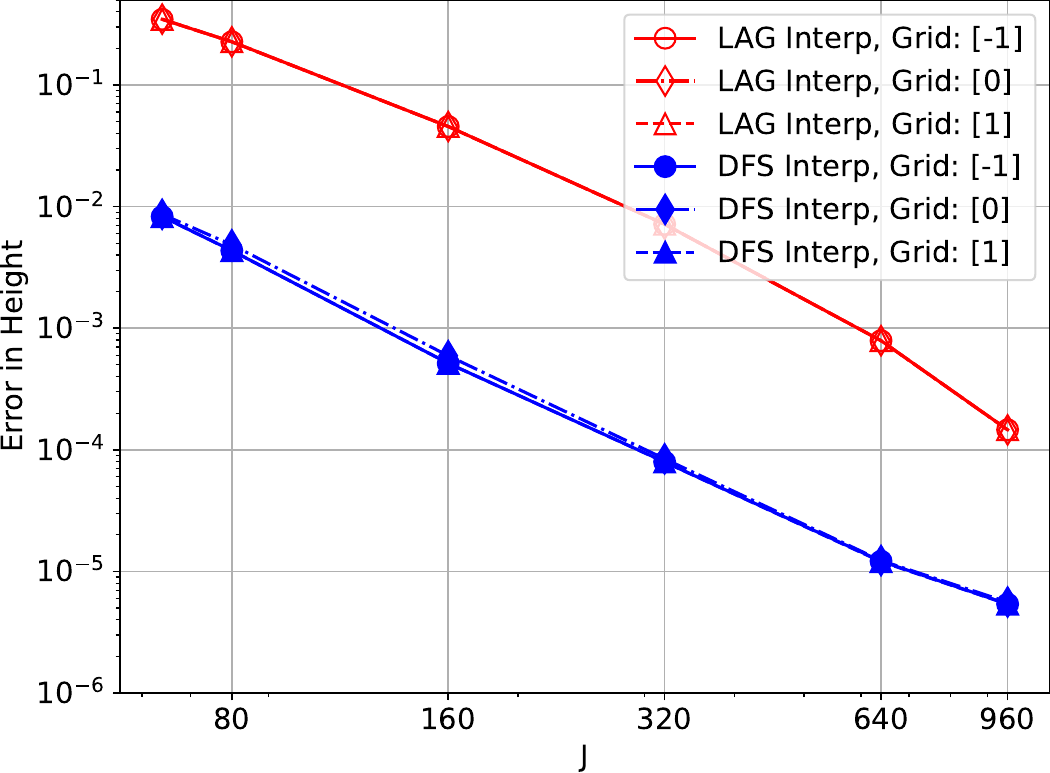}
\caption{$L_1$ norm}
\label{fig:sw1erra}
\end{subfigure}
\hfill
\begin{subfigure}[b]{0.45\textwidth}
\centering
\includegraphics[width=\textwidth]{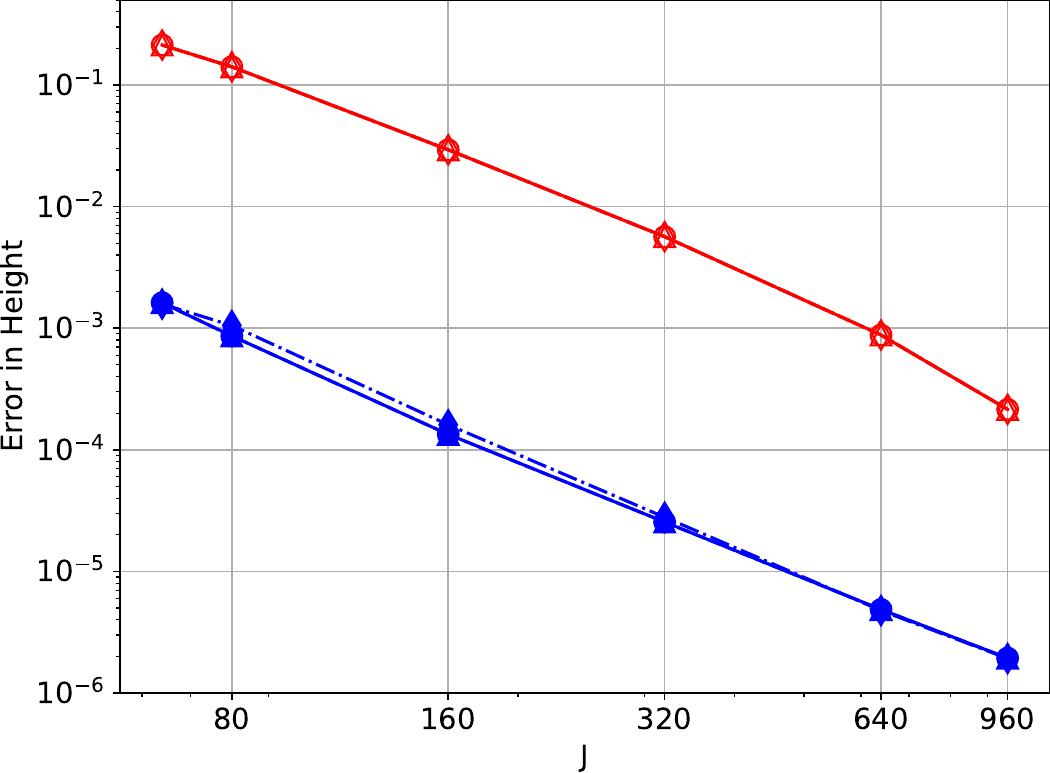}
\caption{$L_2$ norm}
\label{fig:sw1errb}
\end{subfigure}
\hfill \\[0.8em]
\begin{subfigure}[b]{0.45\textwidth}
\centering
\includegraphics[width=\textwidth]{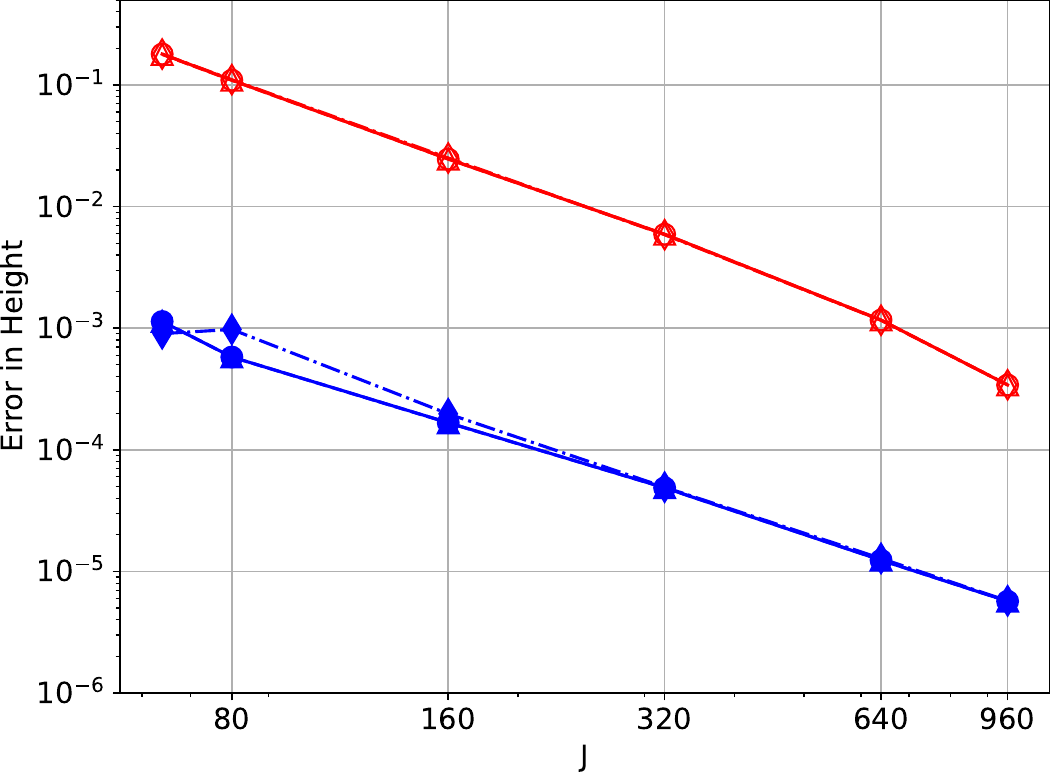}
\caption{Max norm}
\label{fig:sw1errc}
\end{subfigure}
\hfill
\begin{subfigure}[b]{0.45\textwidth}
\centering
\includegraphics[width=\textwidth]{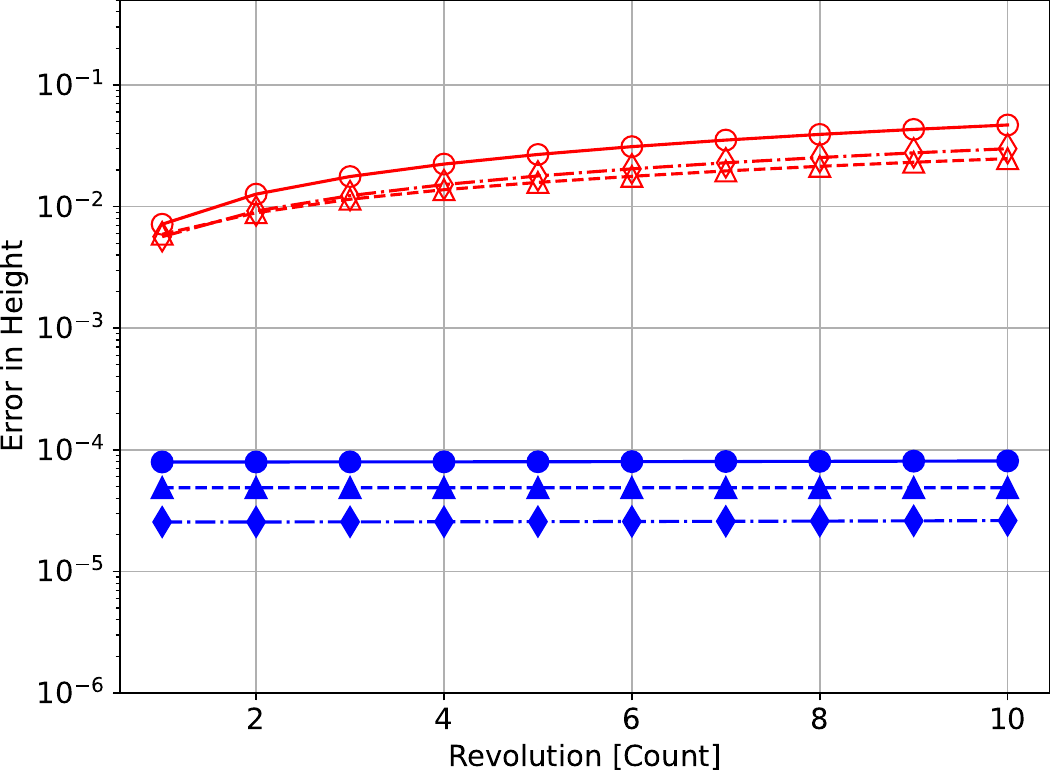}
\caption{$L_2$ error vs number revolutions}
\label{fig:sw1errlong}
\end{subfigure}
\caption{Williamson TC1.  (a)--(c) Relative errors under different norms in the numerical solutions of the predicted height field after one revolution against the latitude grid size $J$. (d) Relative $L_2$ error in the predicted height field as a function of the number of revolutions using Grid [-1] and $J=320$. All experiments were conducted with a time step of 600 s.}
\label{fig:sw1sperr}
\end{figure}
Figures \ref{fig:sw1erra}, \ref{fig:sw1errb}, and \ref{fig:sw1errc} present the relative $L_1$, $L_2$, and max norm errors for the 12-day integration of TC1, comparing the DFS Interp and LAG Interp methods. The results were computed using a truncation wavenumber of $N=J-2$. The plots show that the errors are nearly identical across different grid types for a fixed interpolation method. Notably, we observe that the DFS interpolation method maintains the same spatial convergence rate as the LAG method but yields significantly lower errors (around two orders of magnitude). The convergence rates are algebraic since the cosine bell has only one continuous derivative. 

\begin{figure}[tb]
\centering
\begin{subfigure}[b]{0.45\textwidth}
\centering
\includegraphics[width=\textwidth]{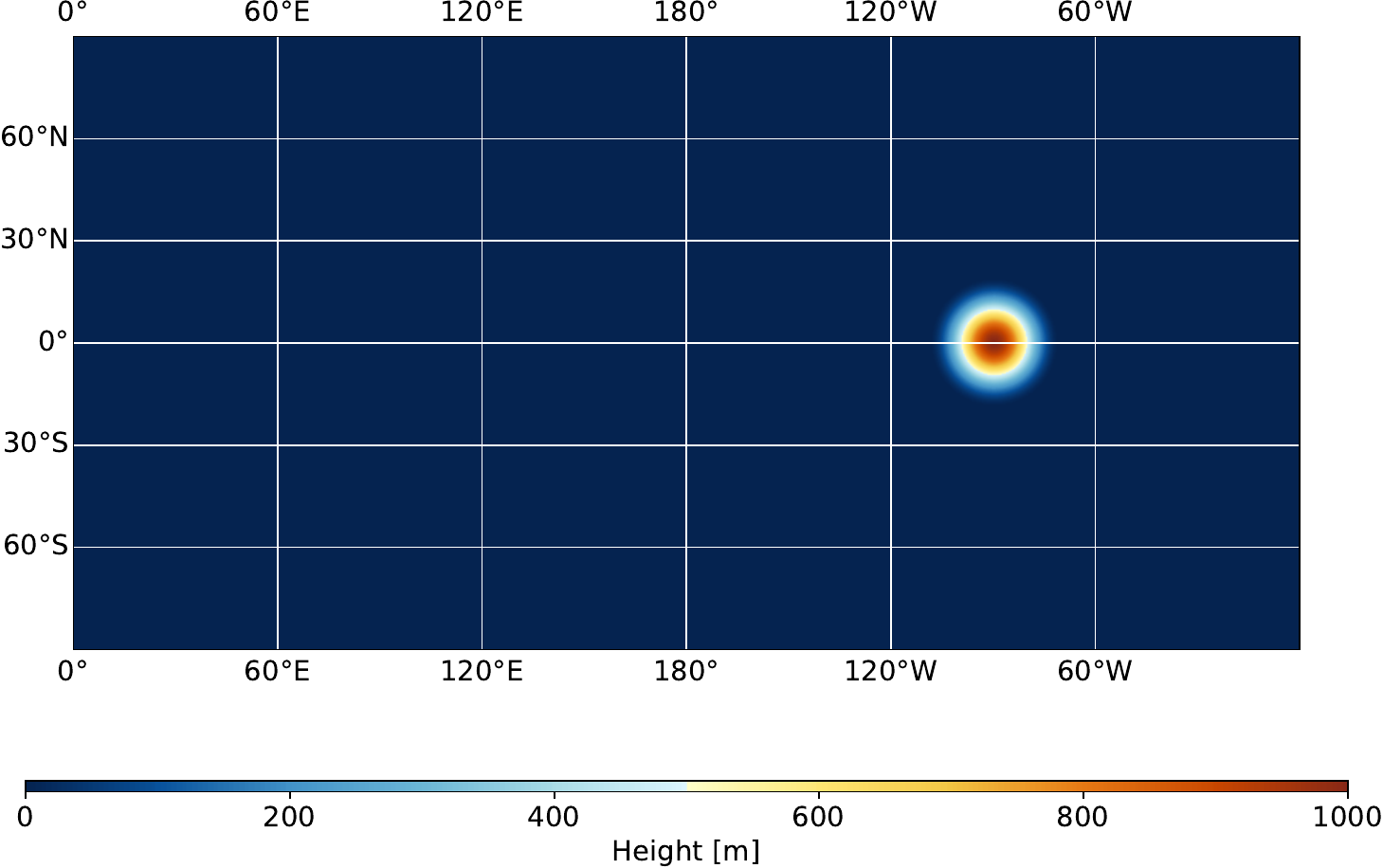}
\caption{DFS Method}
\label{fig:sw1new}
\end{subfigure}
\hfill
\begin{subfigure}[b]{0.45\textwidth}
\centering
\includegraphics[width=\textwidth]{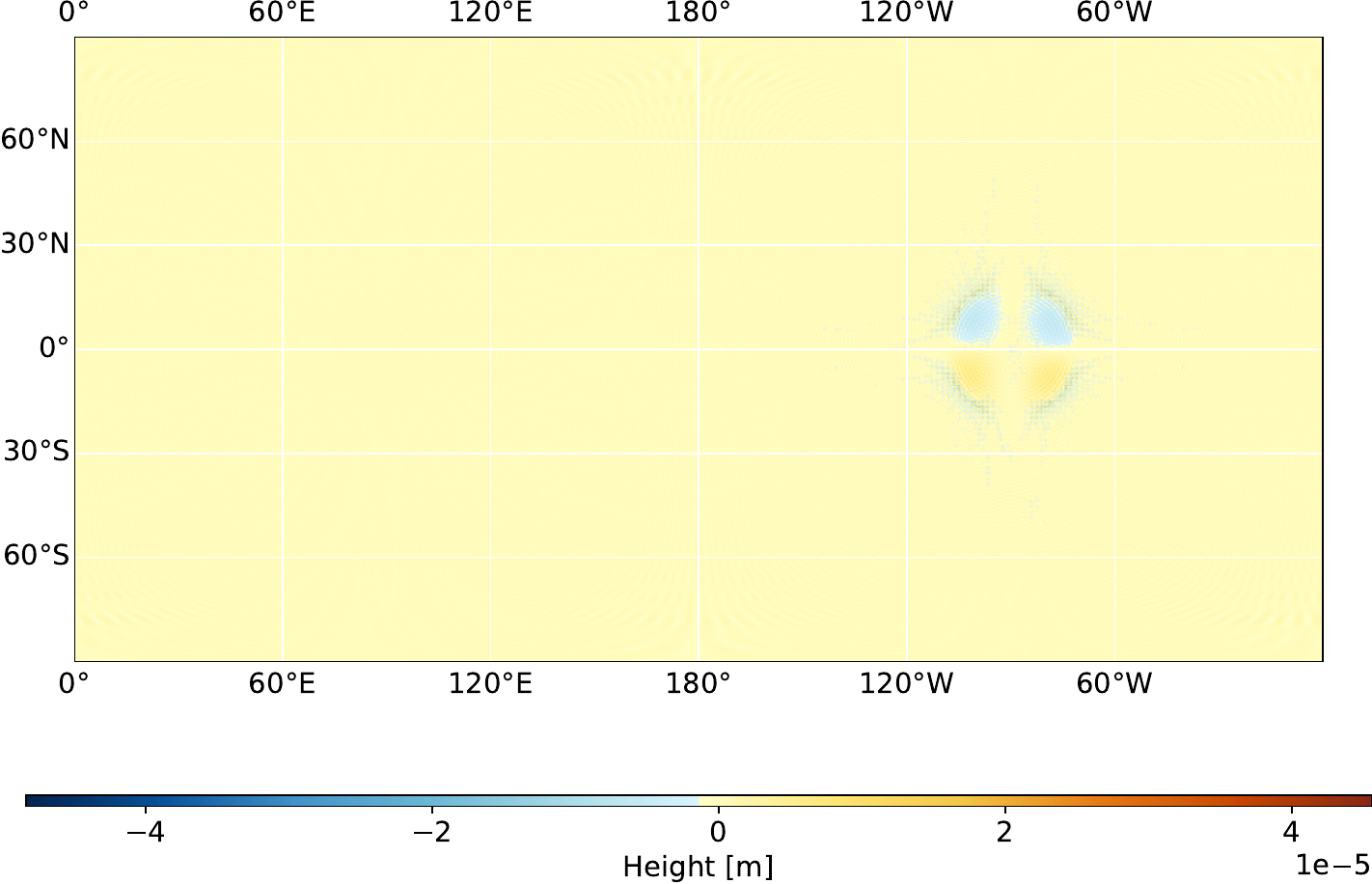}
\caption{DFS Method}
\label{fig:sw1newerr1}
\end{subfigure} 
\hfill \\[0.5em]
\begin{subfigure}[b]{0.45\textwidth}
\centering
\includegraphics[width=\textwidth]{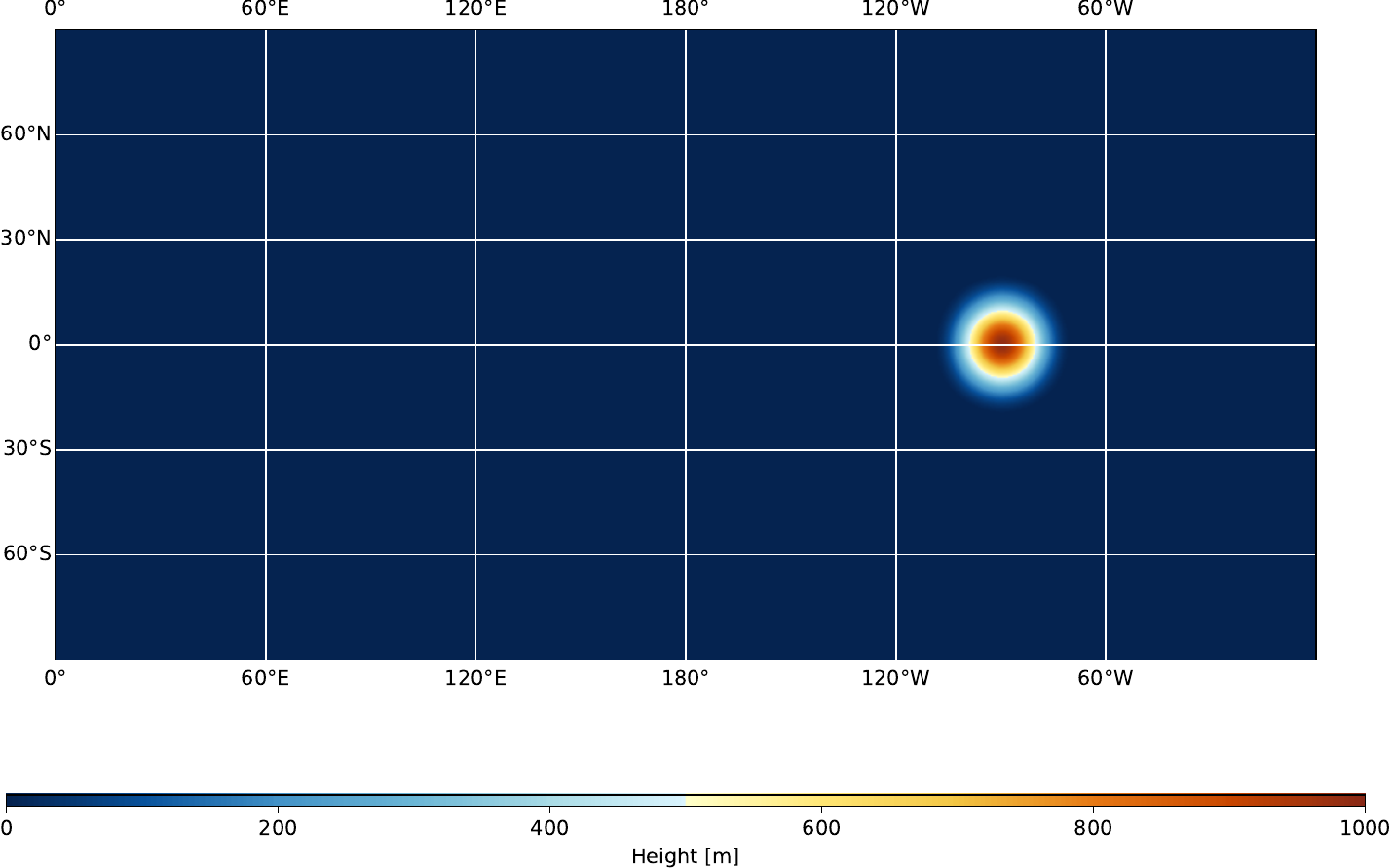}
\caption{LAG Method}
\label{fig:sw1ys}
\end{subfigure}
\hfill
\begin{subfigure}[b]{0.45\textwidth}
\centering
\includegraphics[width=\textwidth]{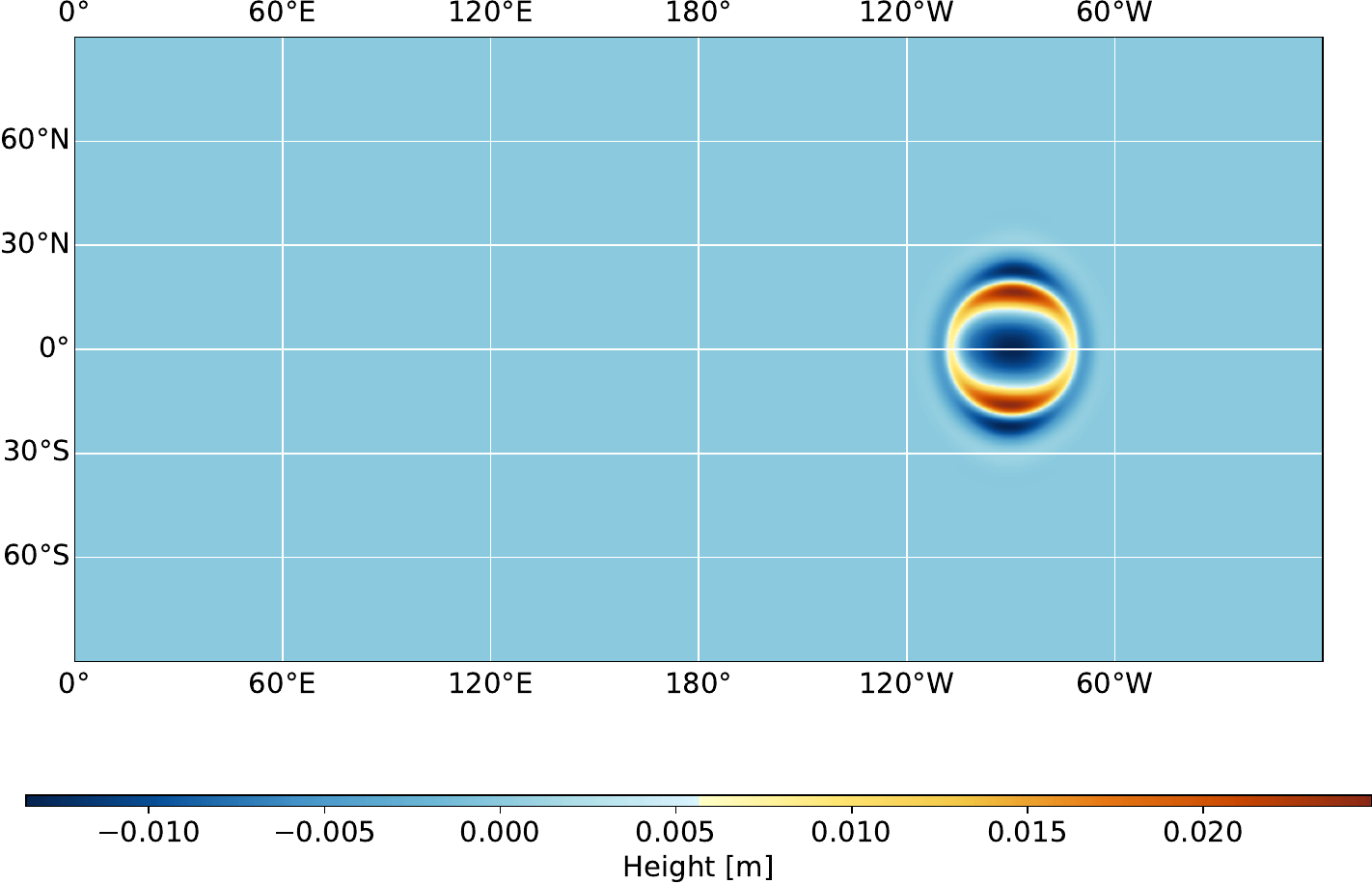}
\caption{LAG Method}
\label{fig:sw1yserr}
\end{subfigure}
\caption{Williamson TC1. (a) \& (c) Heat maps of the predicted height field after 120 days of integration, using the DFS and LAG interpolation methods, respectively. (b) \& (d) Heat maps of the corresponding errors in (a) \& (c), respectively. Experiments were conducted with Grid [-1] using $J=320$ and a time step size of 600 s.}
\label{fig:sw1long}
\end{figure}

To assess the methods under longer integration times, we plot in Figure \ref{fig:sw1errlong} the relative $L_2$ norm errors in the height field from 1 revolution (12 days) to 10 revolutions (120 days).  We also display the final solution after 10 revolutions for the LAG and DFS interpolation methods and the errors in Figure \ref{fig:sw1long}.  These results show that the DFS method is better at preserving accuracy over long periods (about three orders of magnitude smaller errors) and results in a less diffusive method than the LAG method. We also note that the DFS method is better at conserving mass than the LAG method.  For the results in Figure \ref{fig:sw1long}, the absolute error in the mass after 120 days of integration for the DFS method was $1.21\times 10^{-11}$, while it was only $4.34\times 10^{-4}$ for the LAG method.


\subsection{Williamson Test Case 2 (TC2)}
This test case simulates steady-state nonlinear zonal geostrophic flow. In this setup, the angle between the solid body rotation and the polar axis is given again  by $\alpha = \frac{\pi}{2}-0.05$, and the zonal and meridional components of the Coriolis term are expressed as:
	\begin{equation}
		2\mathbf{\Omega}\times \mathbf{r} = 
		2\Omega a\big( \cos \theta\cos\alpha +\cos\phi\sin\theta \sin\alpha, -\sin\phi\sin\alpha
		\big),
	\end{equation}
where $a$ is the Earth's radius.  An exact solution exists for this test case, so the errors for different schemes can be computed over the entire integration time.

\begin{figure}[htb]
\centering
\begin{subfigure}[b]{0.45\textwidth}
\centering
\includegraphics[width=\textwidth]{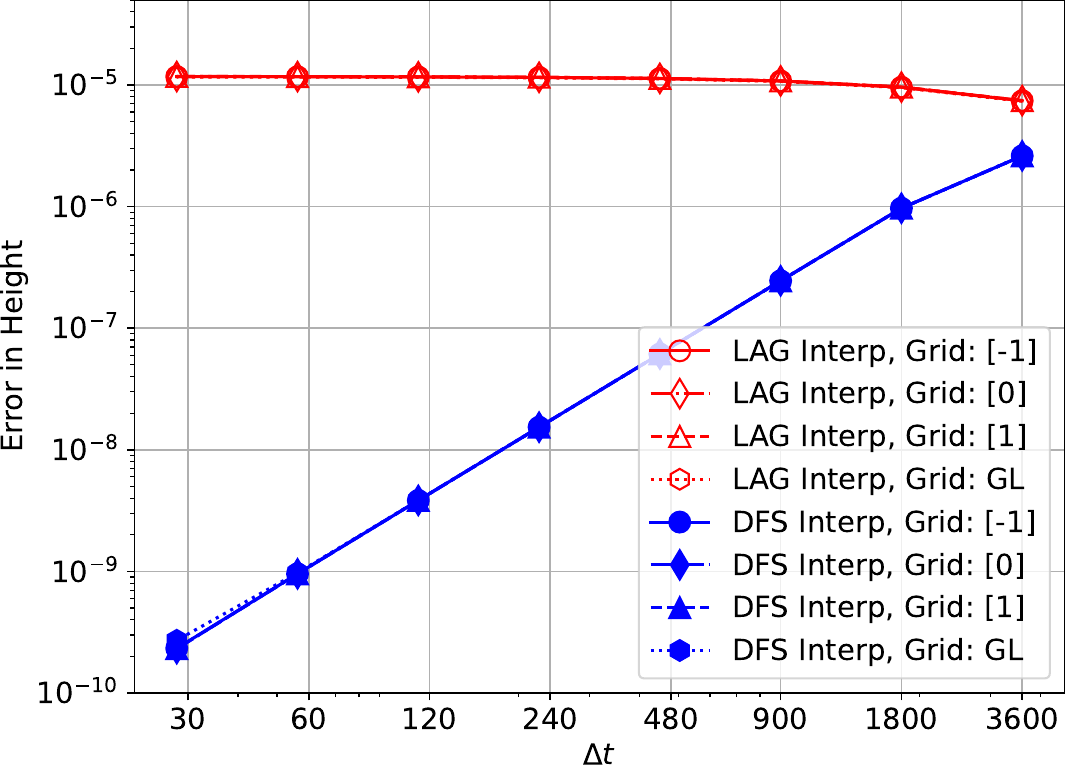}
\caption{$L_1$ norm}
\label{fig:sw2l1}
\end{subfigure}
\hfill
\begin{subfigure}[b]{0.45\textwidth}
\centering
\includegraphics[width=\textwidth]{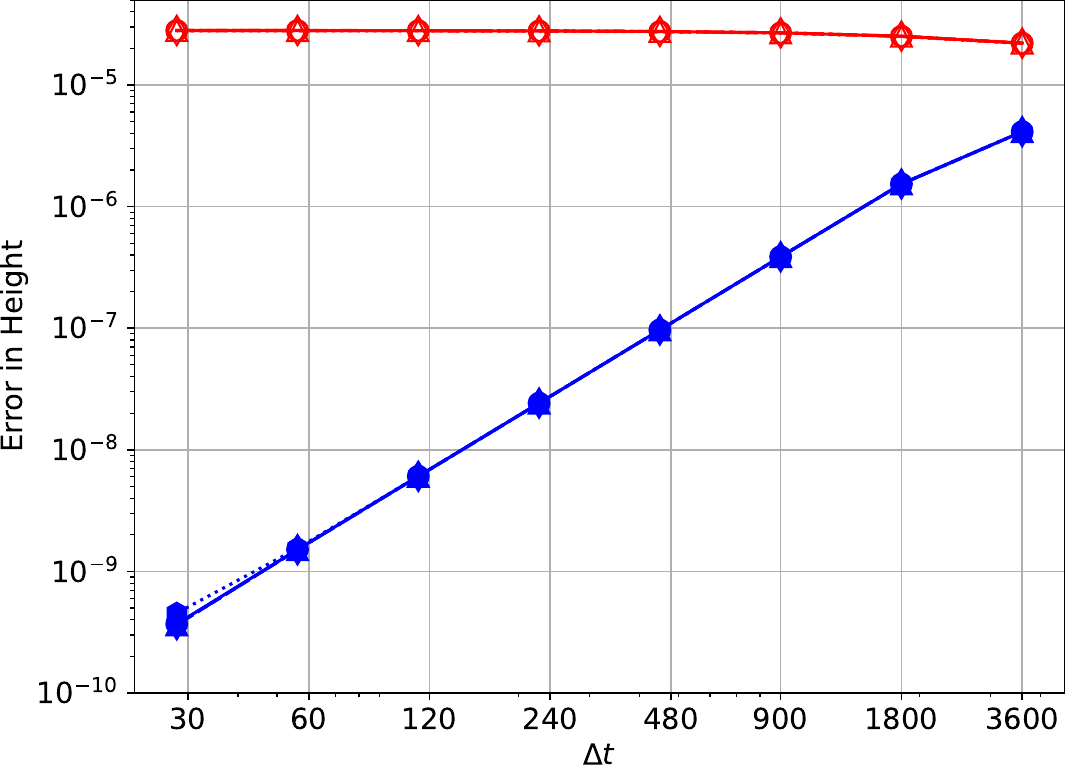}
\caption{Max norm}
\label{fig:sw2max}
\end{subfigure}
\caption{Williamson TC2. Relative $L_1$ and max norm errors for the height field after a 5-day integration of TC2 using $J=80$ and a truncation wavenumber of $N=J-2$. The time step is given in seconds.\label{fig:sw2err}}
\end{figure}

Figure \ref{fig:sw2err} shows the relative $L_1$ and max norms of the errors after a 5-day integration of this test case as the time step is refined.  The $L_2$ errors were between the $L_1$ and max norms, and are thus omitted. 
The DFS SISL and SHT SISL schemes are both expected to be second-order accurate in time. However, we see in Figure \ref{fig:sw2err} that only the schemes that use the DFS interpolation method produce errors that behave at this expected rate.  When using LAG interpolation, neither of the schemes converges in time, which may be due to the spatial errors associated with the scheme dominating the temporal errors.

\begin{figure}[htb]
\centering
\begin{subfigure}[b]{0.45\textwidth}
\centering
\includegraphics[width=\textwidth]{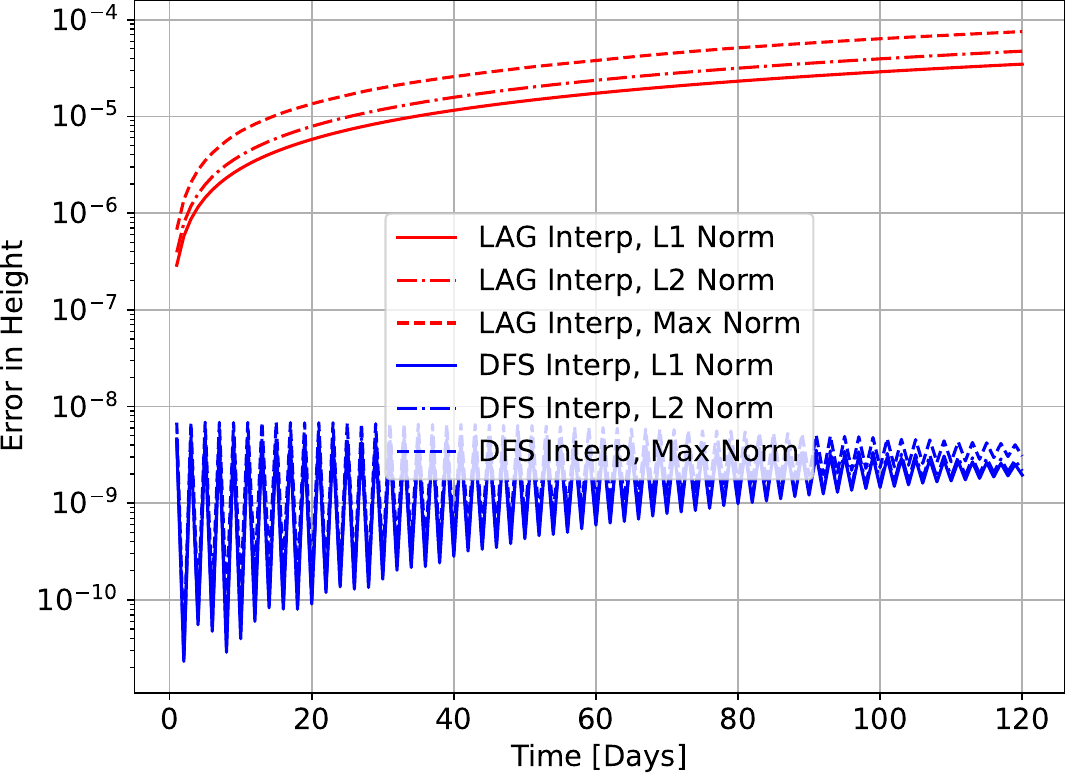}
\caption{Error in Height}
\label{fig:sw2longerror}
\end{subfigure}
\hfill
\begin{subfigure}[b]{0.45\textwidth}
\centering
\includegraphics[width=\textwidth]{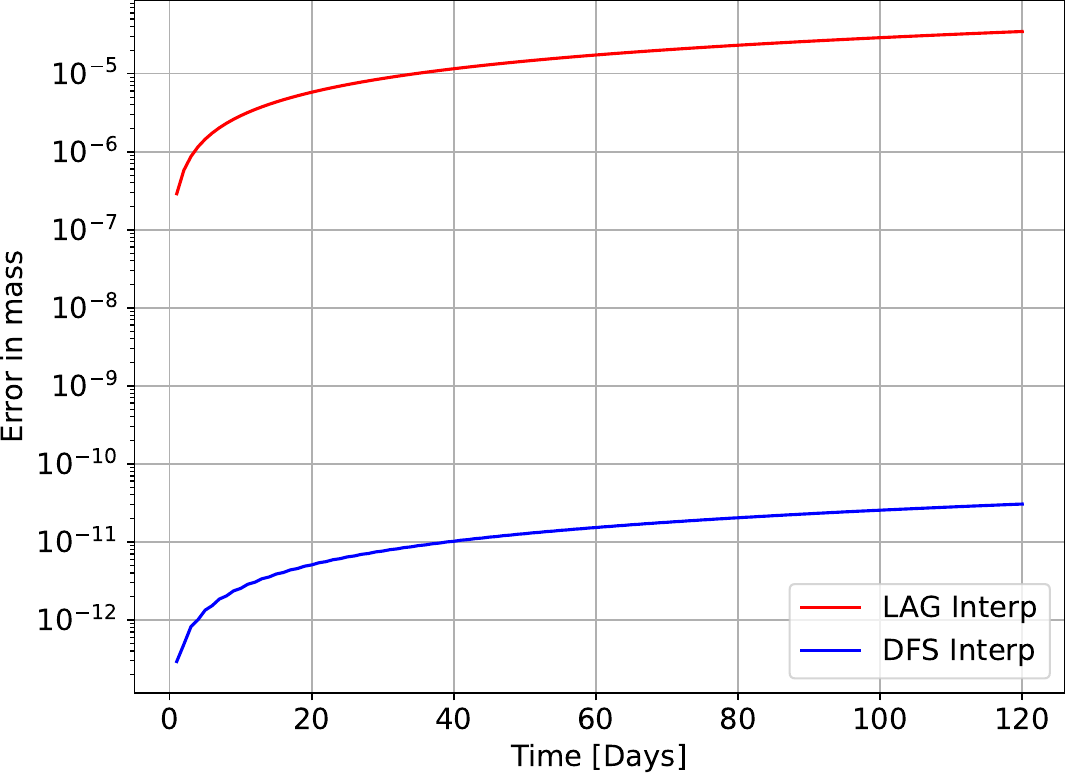}
\caption{Error in Mass}
\label{fig:longmass}
\end{subfigure} 
\caption{Williamson TC2. Time series of errors for the (a) height field and (b) mass over a 120-day integration period using Grid [0] with $J=320$, a truncation wavenumber of $N=J-2$, and a time step of 120 s.}
\label{fig:sw2long}
\end{figure}
	
Figure~\ref{fig:sw2longerror} compares the accuracy of the DFS SISL scheme using the two interpolation methods over a 120-day integration time. We only show results for the DFS SISL scheme on Grid [0] since nearly identical results were obtained for the other grids and for the SHT SISL scheme. Similar to TC1, we see that the time series of the errors in the height field are significantly lower for DFS interpolation.  Additionally, while the errors of the LAG method increase with time, the errors for the DFS method remain relatively stable, with oscillations diminishing over time, indicating that the latter method is better at maintaining geostrophic balance. Figure~\ref{fig:longmass} displays a time series of the errors in the mass from the two interpolation methods.  We again observe that the DFS method demonstrates far superior conservation properties compared to the LAG method. While not displayed, we note that we observed similar results for the errors in energy conservation. 

\subsection{Williamson Test Case 5 (TC5) \label{sec:tc5}}
This test case simulates zonal flow over an isolated conical mountain.  This is a challenging test case as the mountain is impulsively added at $30^{\circ}$N latitude and $90^{\circ}$W longitude to an initially balanced flow, which generates significant gravity waves.  The mountain is specified through the $h_s$ term in \eqref{eq:swe} and is continuous but not differentiable, making convergence tests difficult~\cite{galewsky2004initial}.  This test case does not have an exact solution, so we use two high-resolution solutions computed with the SHT SISL model with a resolution of $J=960$, a truncation wavenumber $N=J-2$, and a time step size of 60 s as reference solutions. One reference solution is computed with LAG interpolation and the other with DFS interpolation.  We interpolate the high-resolution solutions to the lower resolution grids to compute the errors.

\begin{figure}[tb]
\centering
\begin{subfigure}[b]{0.45\textwidth}
\centering
\includegraphics[width=\textwidth]{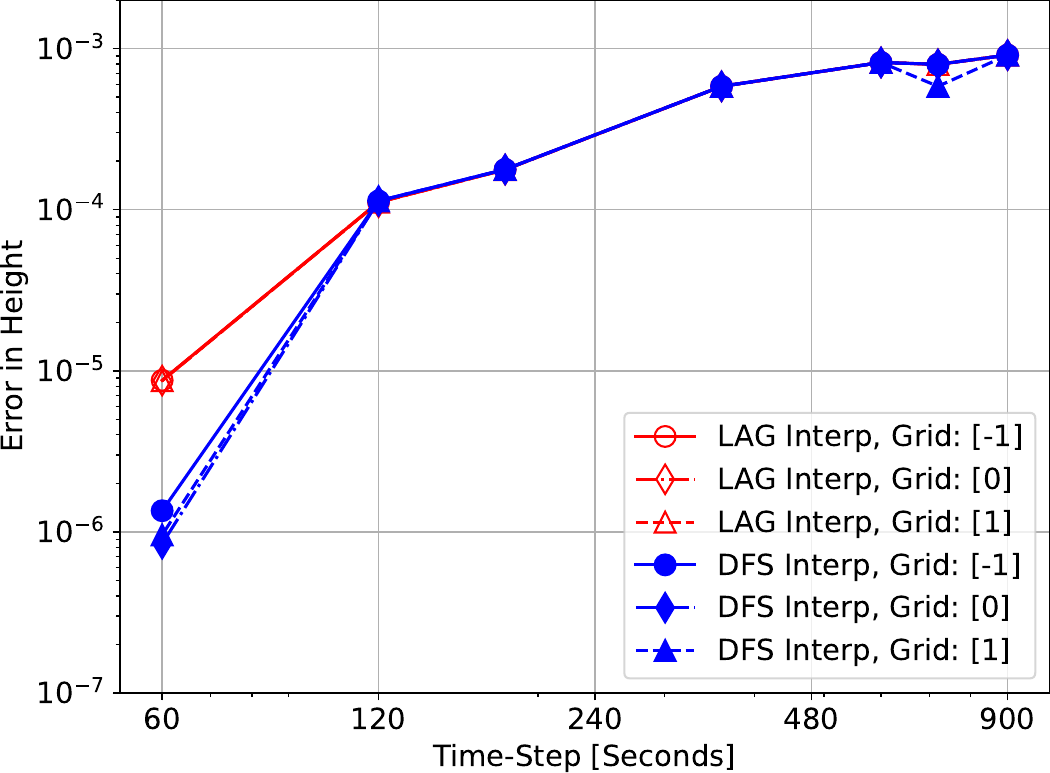}
\caption{LAG interpolation reference }
\label{fig:sw5tlag}
\end{subfigure}
\hfill
\begin{subfigure}[b]{0.45\textwidth}
\centering
\includegraphics[width=\textwidth]{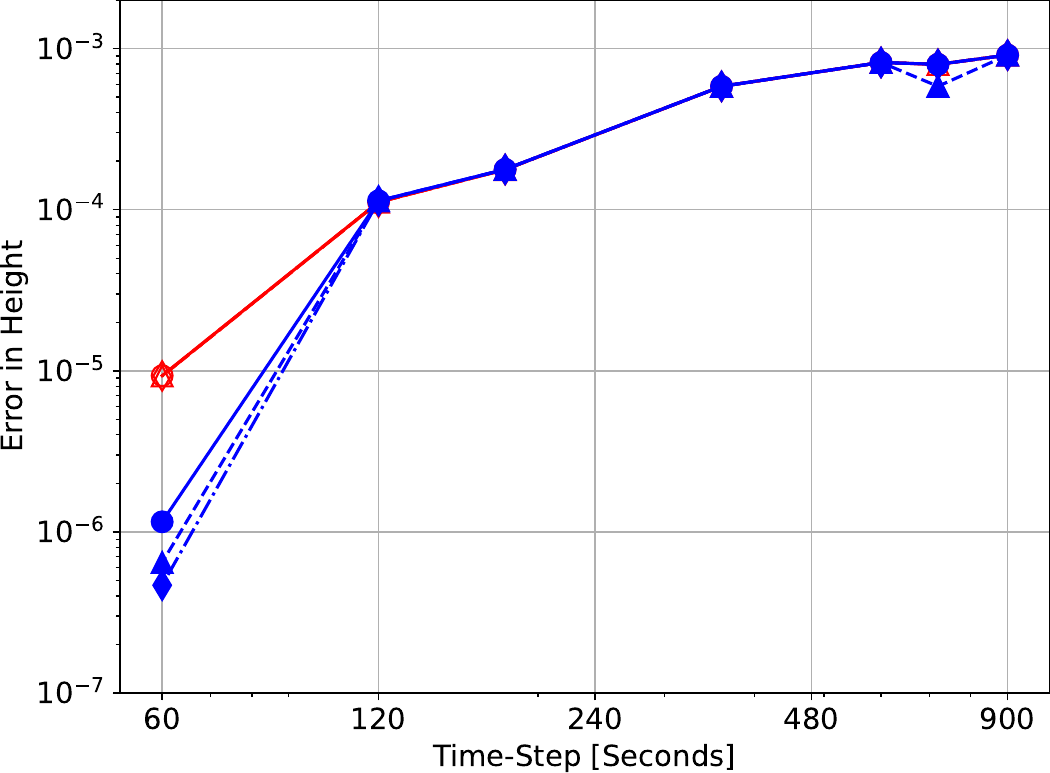}
\caption{DFS interpolation reference }
\label{fig:sw5tdfs}
\end{subfigure} 
\caption{Williamson TC5. Relative $L_2$ norm errors for the height field at 15 days of integration as the time step is decreased.  All results are for a grid size with $J=320$ truncation $N=J-2$.  Errors computed using high-resolution SHT SISL solutions with (a) LAG and (b) DFS interpolation.}
\label{fig:sw5err}
\end{figure}

Figure~\ref{fig:sw5err} shows the $L_2$ norm of the errors in the height field after 15 days of integration as the time step is refined.  The solutions are computed with a grid size of $J=320$ and the same $N=J-2$ truncation wavenumber as the reference solution . Errors in Figure~\ref{fig:sw5tlag} were computed using the LAG interpolation reference, while errors in Figure~\ref{fig:sw5tdfs} were computed using the DFS interpolation reference. Comparing both reference solutions, we observe that the errors are nearly identical for the DFS and LAG methods when the time step size exceeds 120 s, regardless of the grid type used. However, the error drops considerably for $\Delta t=60$ s, with the DFS method producing smaller errors compared to the LAG method. This error drop is most likely due to the fact that the reference solutions were computed at $\Delta t=60$ s, so the computed solutions are capturing the same gravity waves generated by the initial imbalance from the mountain.

\begin{figure}[t]
\centering
\begin{subfigure}[b]{0.45\textwidth}
\centering
\includegraphics[width=\textwidth]{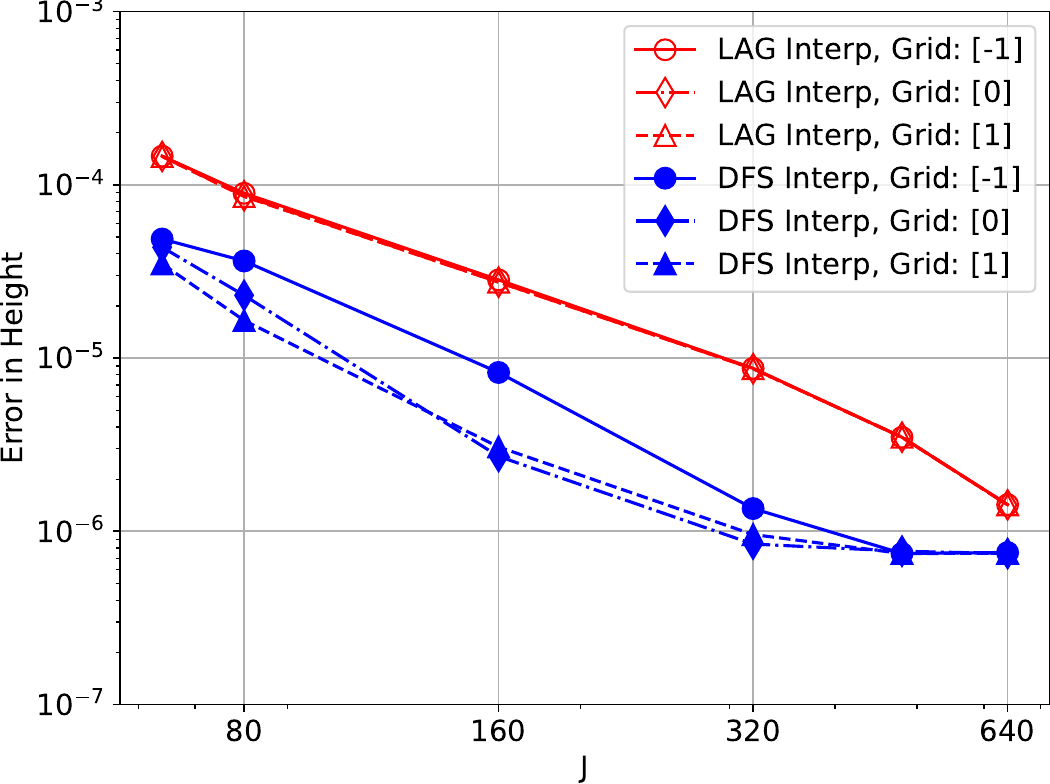}
\caption{LAG interpolation reference}
\label{fig:sw5splag}
\end{subfigure}
\hfill
\begin{subfigure}[b]{0.45\textwidth}
\centering
\includegraphics[width=\textwidth]{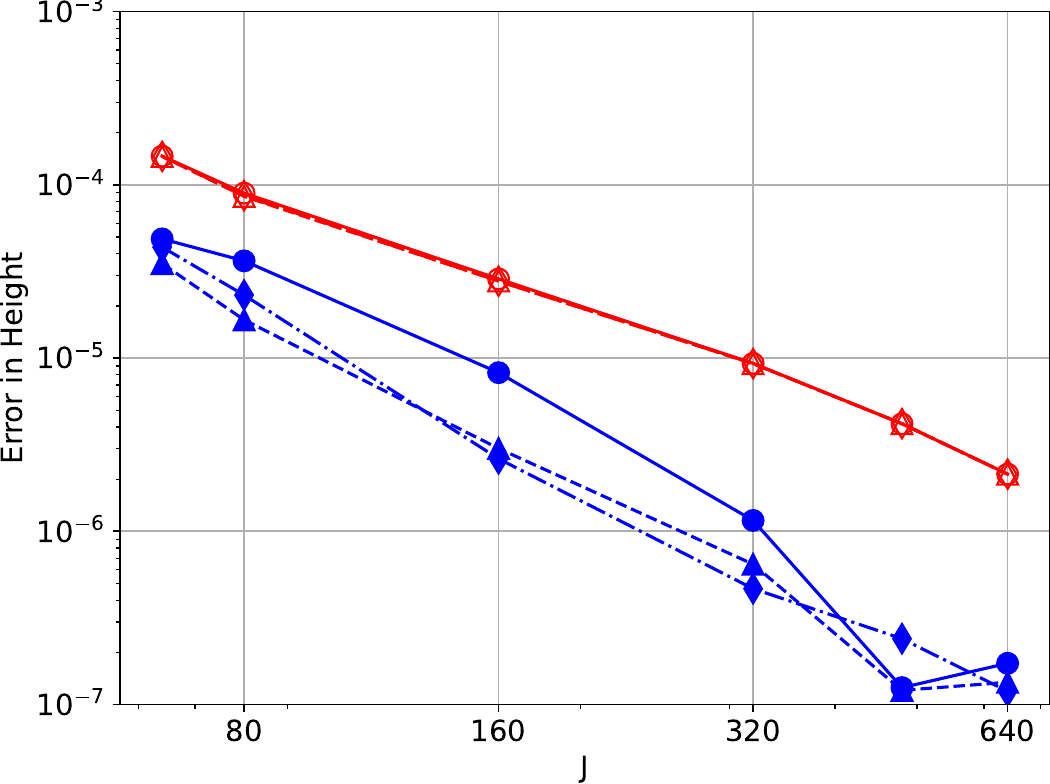}
\caption{DFS interpolation reference}
\label{fig:sw5spdfs}
\end{subfigure}
\caption{Williamson TC5. Relative $L_2$ errors for the height field after 15 days of integration with a time step of 60 s and increasing spatial resolutions. Errors computed using high-resolution SHT SISL solutions with (a) LAG  and (b) DFS interpolation.}
\label{fig:sw5sp}
\end{figure}

We next examine the spatial convergence of the methods in  Figure \ref{fig:sw5sp} after 15 days of integration. For this test, the time step is set to 60 seconds and $J=64, 80, 160, 320, $ and $640$. For both reference solutions, we see that the DFS method gives lower errors than the LAG method across all grid types and sizes. Comparing the results across the reference solutions, we see that the errors are similar up to $J=160$, but for larger $J$, convergence  stalls for the DFS method when using the LAG reference solution.  However, when switching to the DFS reference, the DFS interpolation schemes continue to converge.  In Figure \ref{fig:sw5spdfs} we see that for Grids [0] and [1], the $L_2$ errors in the height field for the DFS interpolation method are consistently an order of magnitude smaller than the LAG interpolation for $J=160$ to $J=640$.
\begin{figure}[!h]
\begin{subfigure}[b]{0.47\textwidth}
\includegraphics[width=\textwidth]{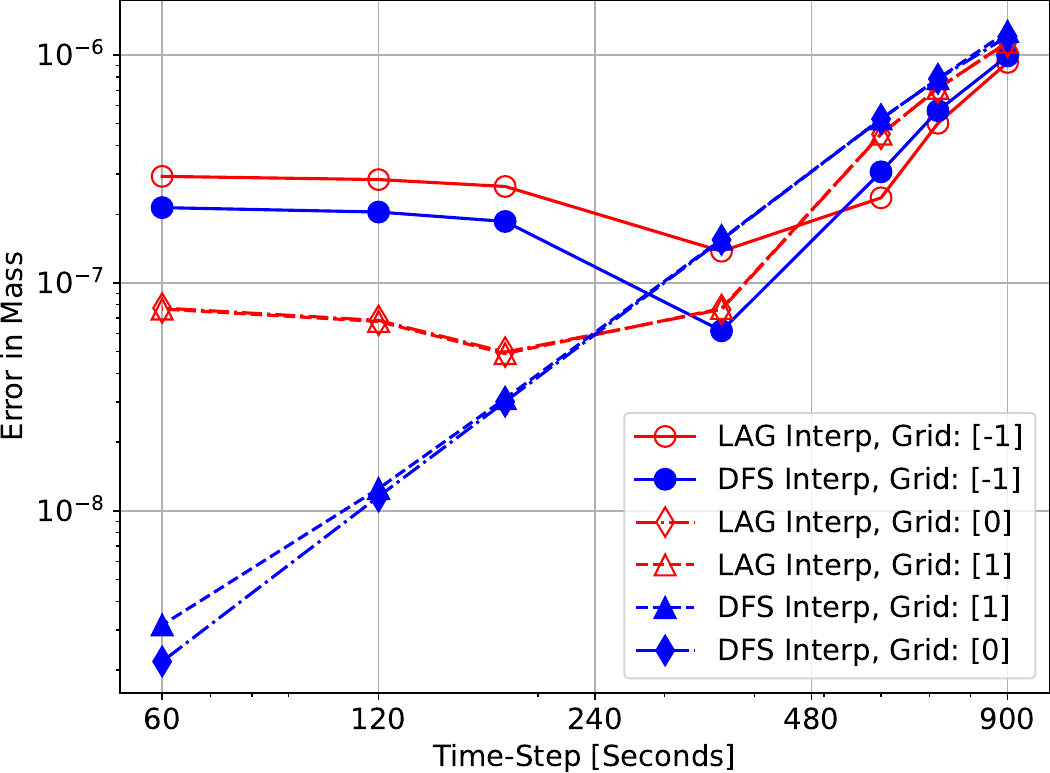}
\caption{Temporal convergence}
\label{fig:sw5masst}
\end{subfigure}
\hfill
\begin{subfigure}[b]{0.47\textwidth}
\includegraphics[width=\textwidth]{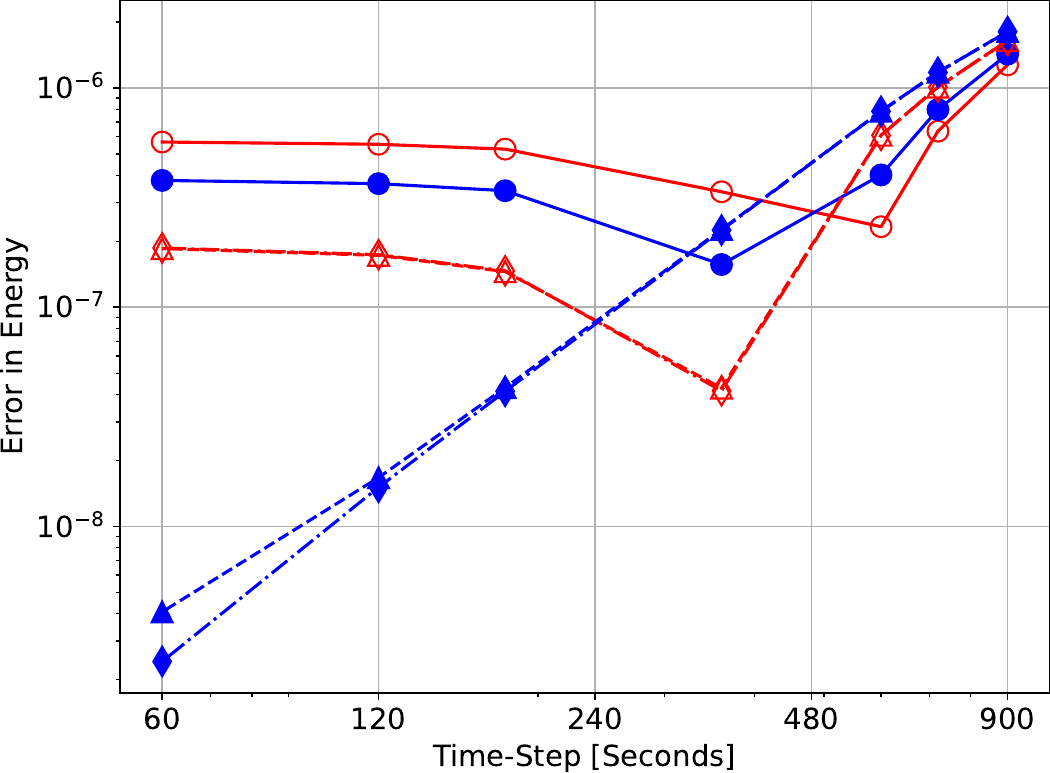}
\caption{Temporal convergence}
\label{fig:sw5energyt}
\end{subfigure}
\hfill \\[0.8em]
\begin{subfigure}[b]{0.47\textwidth}
\includegraphics[width=\textwidth]{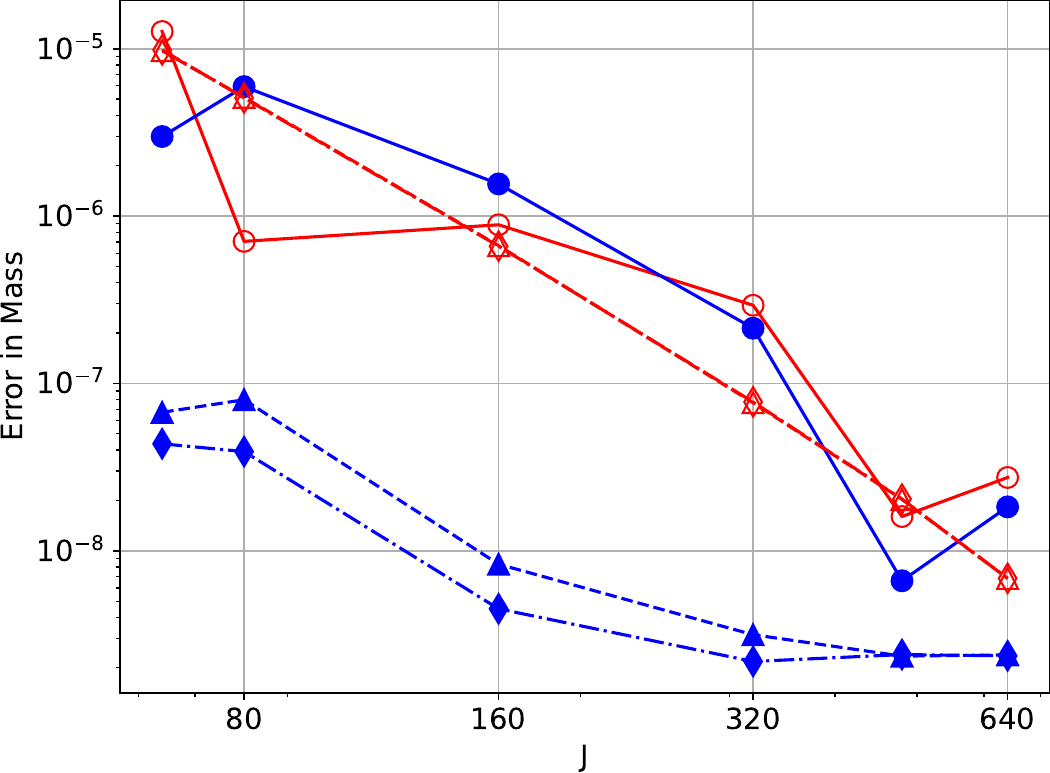}
\caption{Spatial convergence}
\label{fig:sw5massnlon}
\end{subfigure} 
\hfill
\begin{subfigure}[b]{0.47\textwidth}
\includegraphics[width=\textwidth]{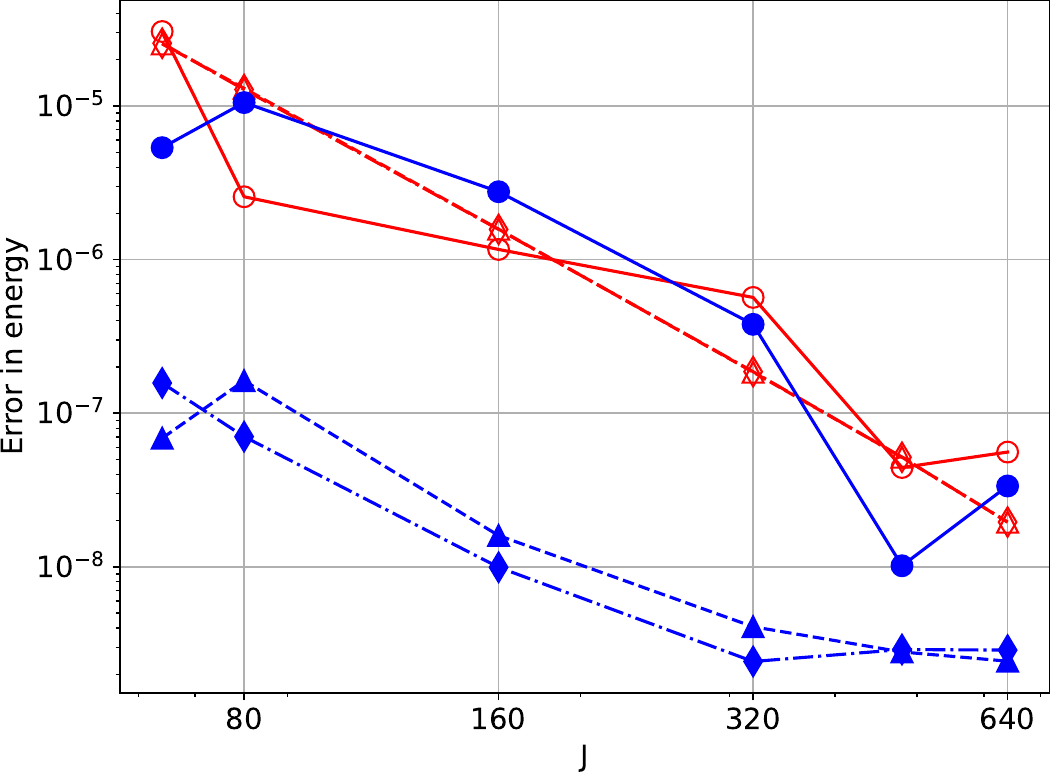}
\caption{Spatial convergence}
\label{fig:sw5energynlon}
\end{subfigure}
\caption{Williamson TC5. (a) Mass and (b) energy conservation errors at day 15 with a fixed spatial resolution of $J=320$ and time steps ranging from 900 s to 60 s. (c) Mass and (d) energy conservation errors at day 15 with a fixed time step of 60 s and spatial resolutions ranging from  $J=64$ to $J=640$.}
\label{fig:mass}
\end{figure}

Convergence results for the errors in mass and energy conservation after 15 days of integration are shown in Figure \ref{fig:mass}, with temporal convergence shown in (a)--(b) and spatial convergence shown in (c)--(d).  For the temporal convergence, we see that initially the LAG and DFS interpolation methods converge at similar rates, but as the time step further decreases, only the DFS method with Grids [1] and [0] continues to converge.  For the spatial convergence, we see that for Grids [0] and [1] the conservation of mass and energy for the DFS interpolation method is again significantly better than that of the LAG interpolation method.  Interestingly, for Grid [-1], the conservation errors for both the LAG and DFS interpolation methods vary widely, unlike those for Grids [0] and [1]. 
\clearpage

\begin{figure}[ht]
\centering
\begin{subfigure}[b]{0.49\textwidth}
\includegraphics[width=\textwidth]{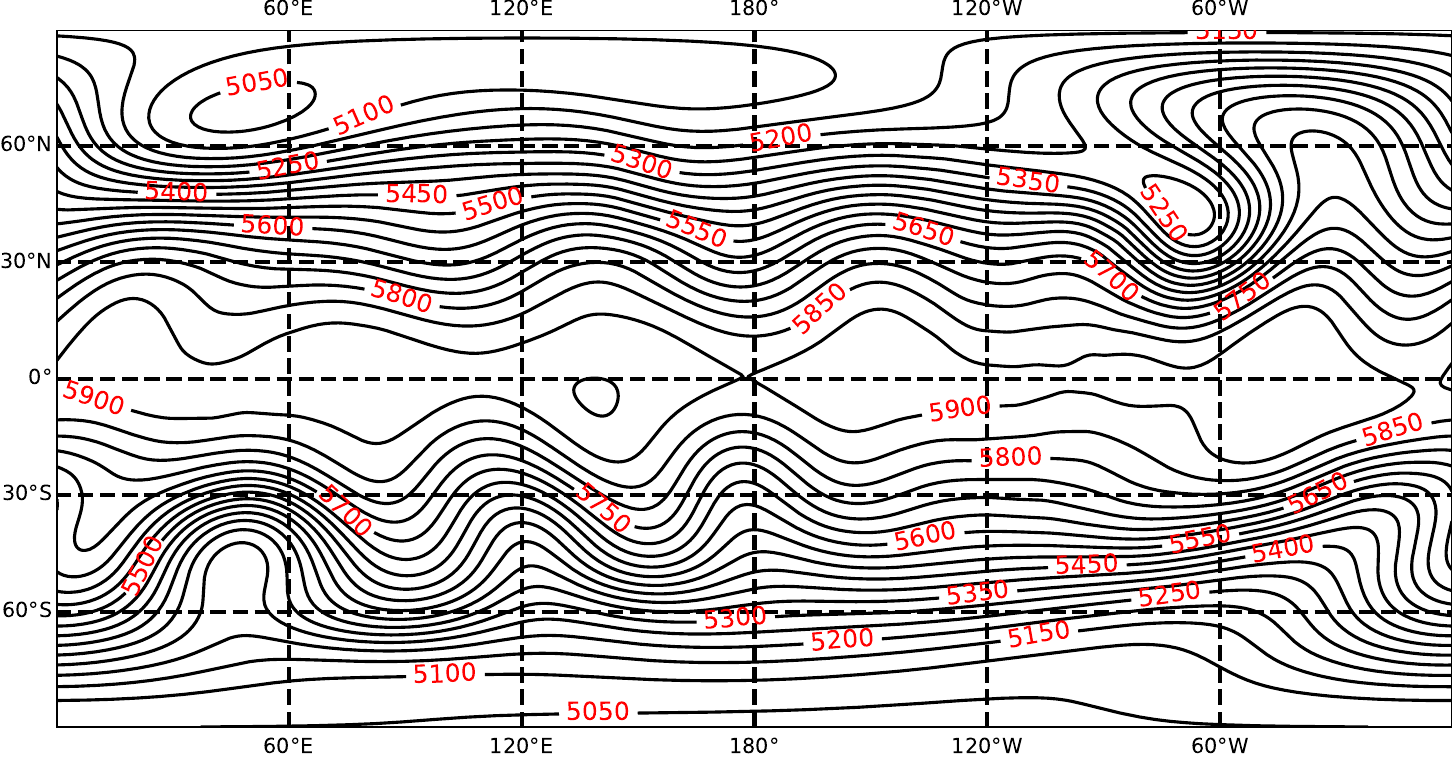}
\caption{Height field, LAG Interp}
\label{fig:sw5lagplot}
\end{subfigure}
\hfill
\begin{subfigure}[b]{0.49\textwidth}
\includegraphics[width=\textwidth]{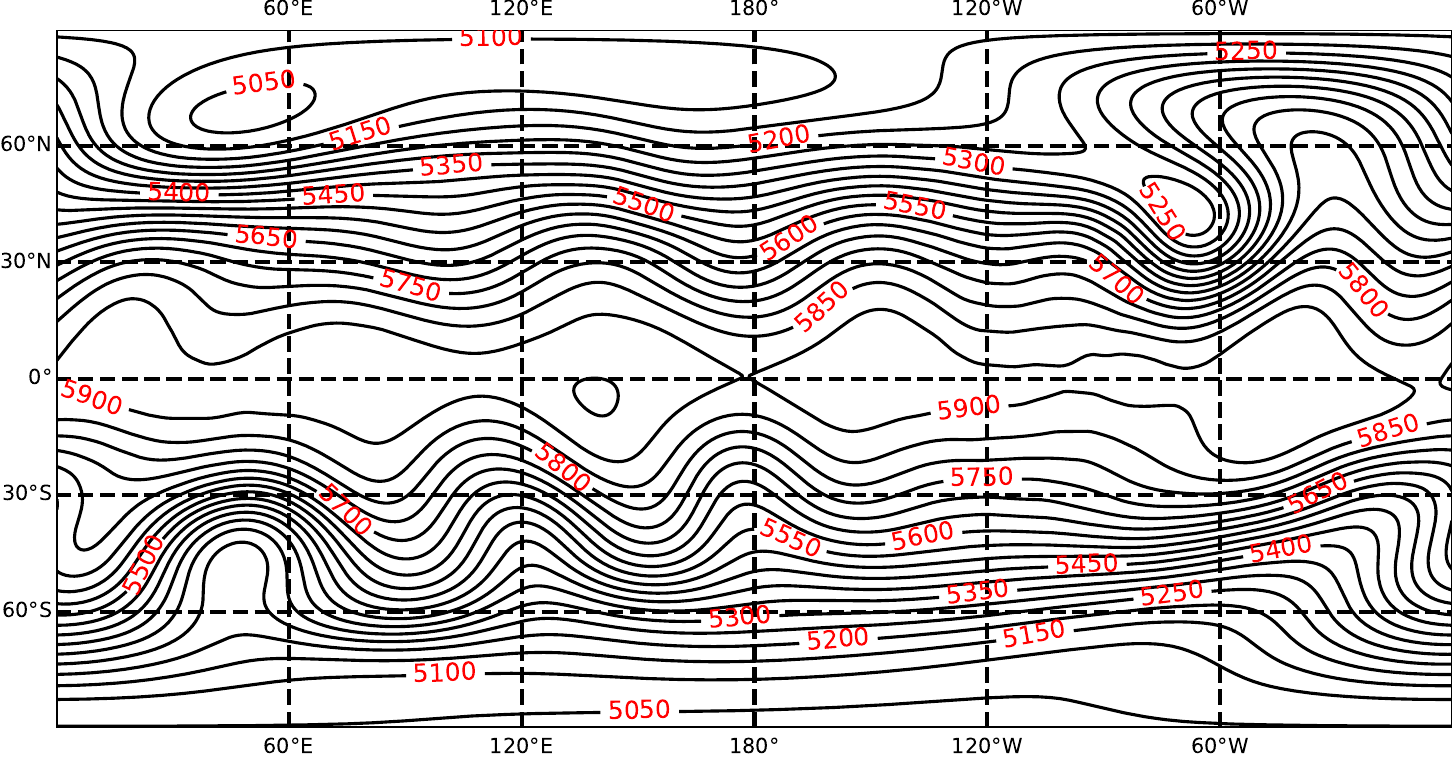}
\caption{Height field, DFS Interp}
\label{fig:sw5dfsplot}
\end{subfigure} \\[0.5em]
\begin{subfigure}[b]{0.49\textwidth}
\includegraphics[width=\textwidth]{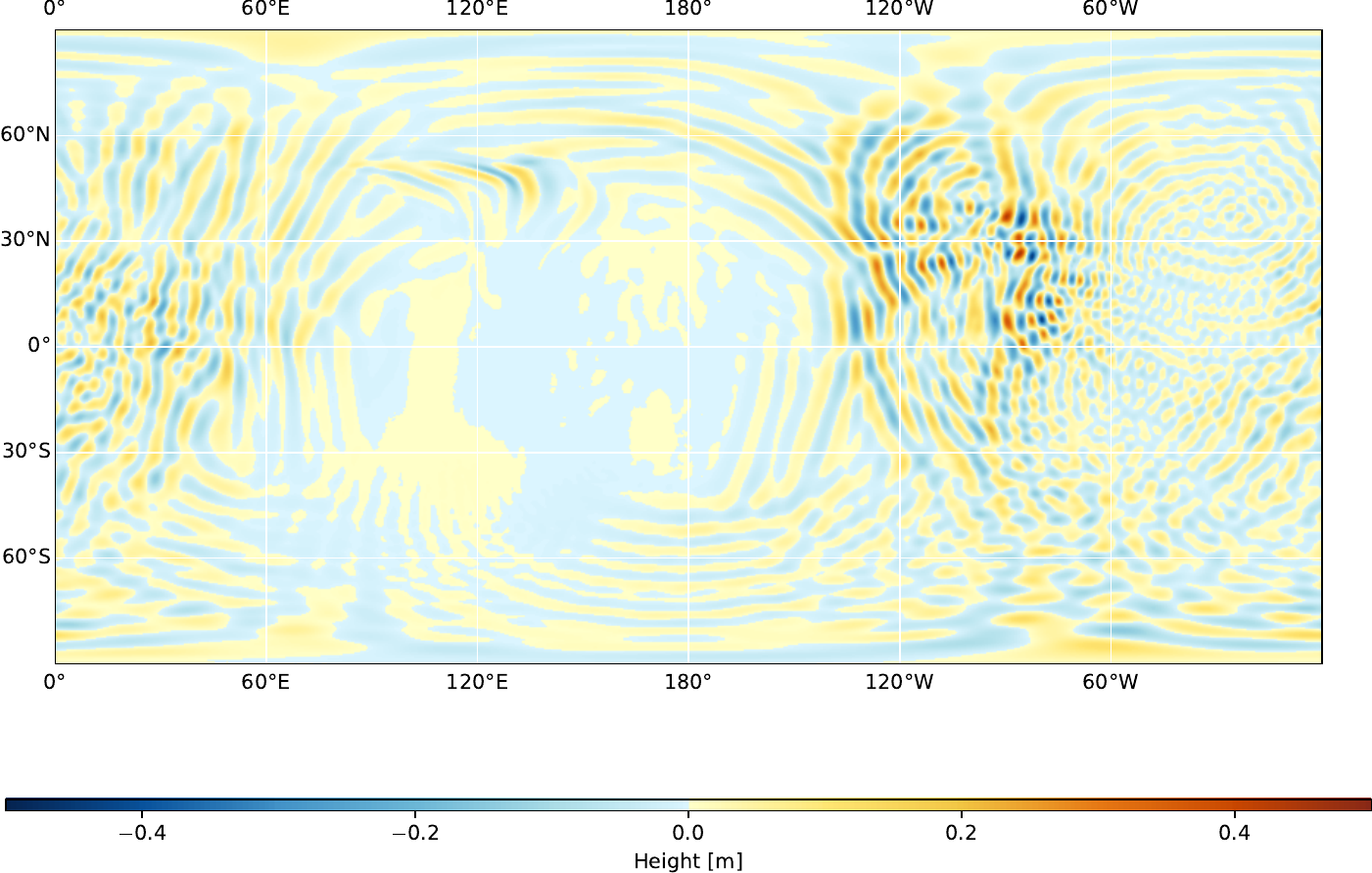}
\caption{Error, LAG Interp}
\label{fig:sw5lagerrplot}
\end{subfigure}
\hfill
\begin{subfigure}[b]{0.49\textwidth}
\includegraphics[width=\textwidth]{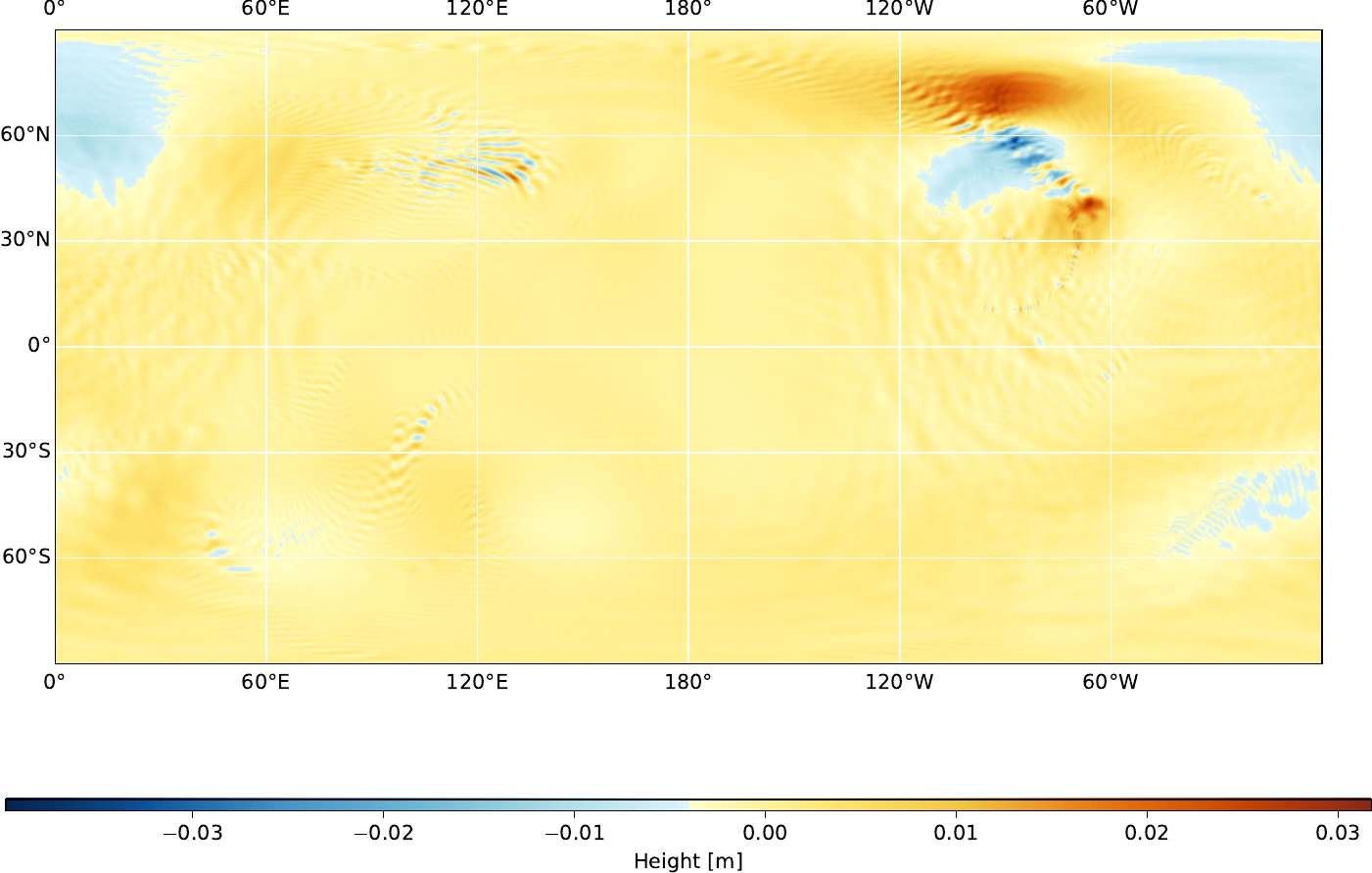}
\caption{Error, DFS Interp}
\label{fig:sw5dfserrplot}
\end{subfigure} 
\caption{Williamson TC5. Comparison of the the computed height field from the two methods using grid [0], $J=320$, and time step of 60 s.  Contours in (a) and (b) are plotted in increments of 50 m starting at $h=5000$m and ending at $h=6000$m.  Errors in (c) and (d) are computed using the high-resolution solutions based on the corresponding interpolation method; note the different color bars.
\label{fig:sw5plot}}
\end{figure}

Figure \ref{fig:sw5plot} (a) \& (b) display contour plots of the height field computed using LAG and DFS interpolation methods, respectively, while (c) \& (d) show the error in these solutions. For the LAG (DFS) results, the errors are measured against the high-resolution SHT solution computed using LAG (DFS) interpolation, respectively.  The height fields look nearly identical, but we see that the error for DFS is an order of magnitude smaller than LAG and is also more localized to the location of the conical mountain.

\subsection{Williamson Test Case 6 (TC6)\label{sec:tc6}}
 This test case is initialized with a wavenumber four Rossby–Haurwitz wave. For the barotropic vorticity equation, this initial condition propagates eastward across the sphere without changing its shape. However, for the SWE, this motion is only approximated, leading to inherent instabilities in the flow~\cite{thuburn2000numerical}.  This test also does not have an exact solution, so similar to TC5, we use high-resolution SHT-SISL schemes to compute the reference solutions using both the LAG and DFS interpolation methods. 

\begin{figure}[htb]
		\begin{subfigure}[b]{0.47\textwidth}
        \includegraphics[width=\textwidth]{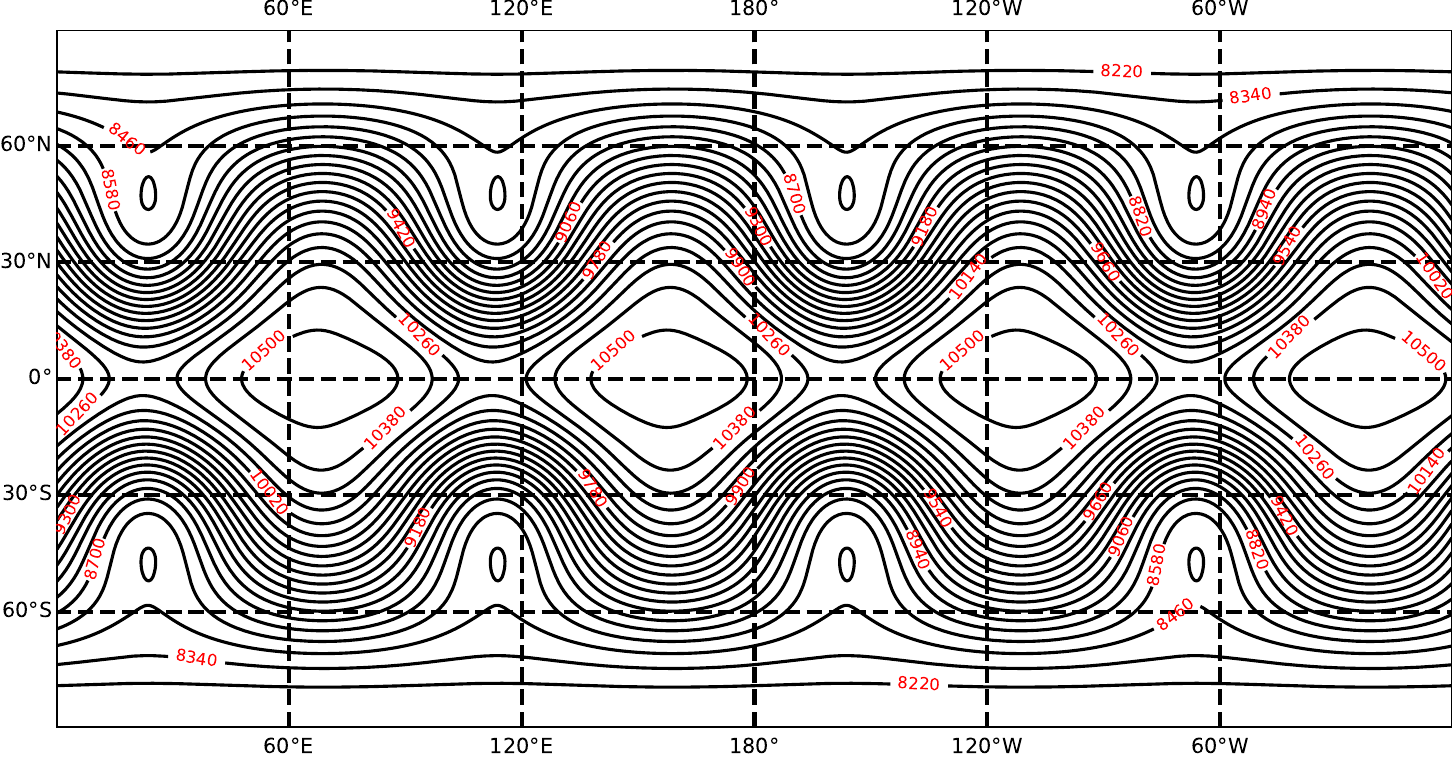}
			\caption{LAG interpolation}
			\label{fig:sw6lagdfsmodel}
		\end{subfigure}
		\hfill
		\begin{subfigure}[b]{0.47\textwidth}
			\includegraphics[width=\textwidth]{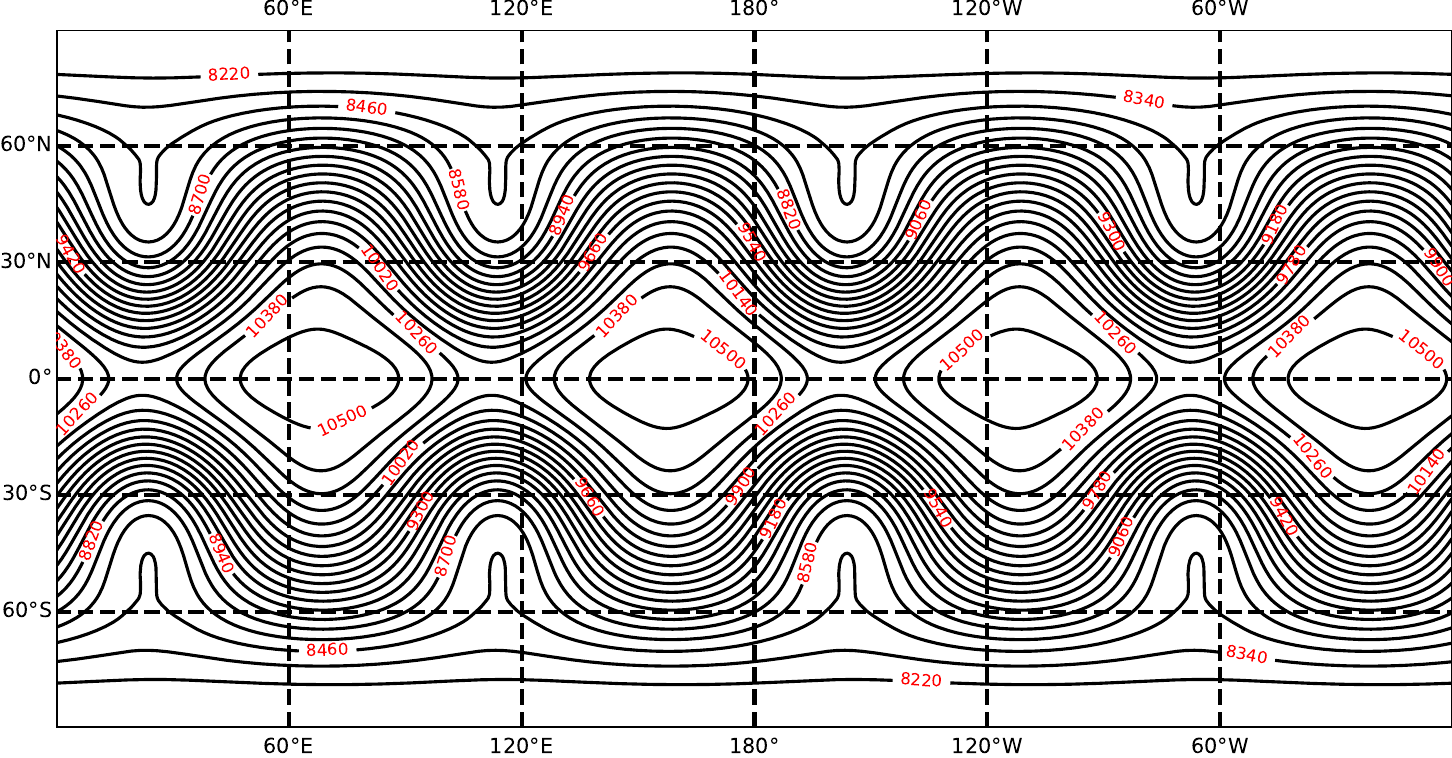}
			\caption{DFS interpolation}
			\label{fig:sw6dfsdfsmodel}
		\end{subfigure}
		\caption{Williamson TC6. Contour plots of the predicted height field after 14 days of integration using the DFS-SISL scheme with Grid [0], $J=320$, and $\Delta t = 60$ s.  Contours range from $8100$m to $10500$m in increments of $120$m. Note the differences in the troughs for the 8460 m contour.}
  \label{fig:sw6contour}
\end{figure}

We first show a qualitative comparison of the solutions computed, which is common for this test case.  Figure~\ref{fig:sw6lagdfsmodel} and \ref{fig:sw6dfsdfsmodel} show contour plots of the height field after 14 days of integration for the Lagrange and DFS interpolation methods, respectively, using the DFS-SISL scheme. We see that the results look similar, except that the $8460$ m contour for the LAG interpolation method is closed in each trough, while it is open for the DFS method.  This latter behavior  better matches the results from~\cite{thuburn2000numerical}, which does a detailed study of TC6 using a high-order Eulerian SHT model.


\begin{figure}[htb]
\centering
\begin{subfigure}[b]{0.45\textwidth}
\includegraphics[width=\textwidth]{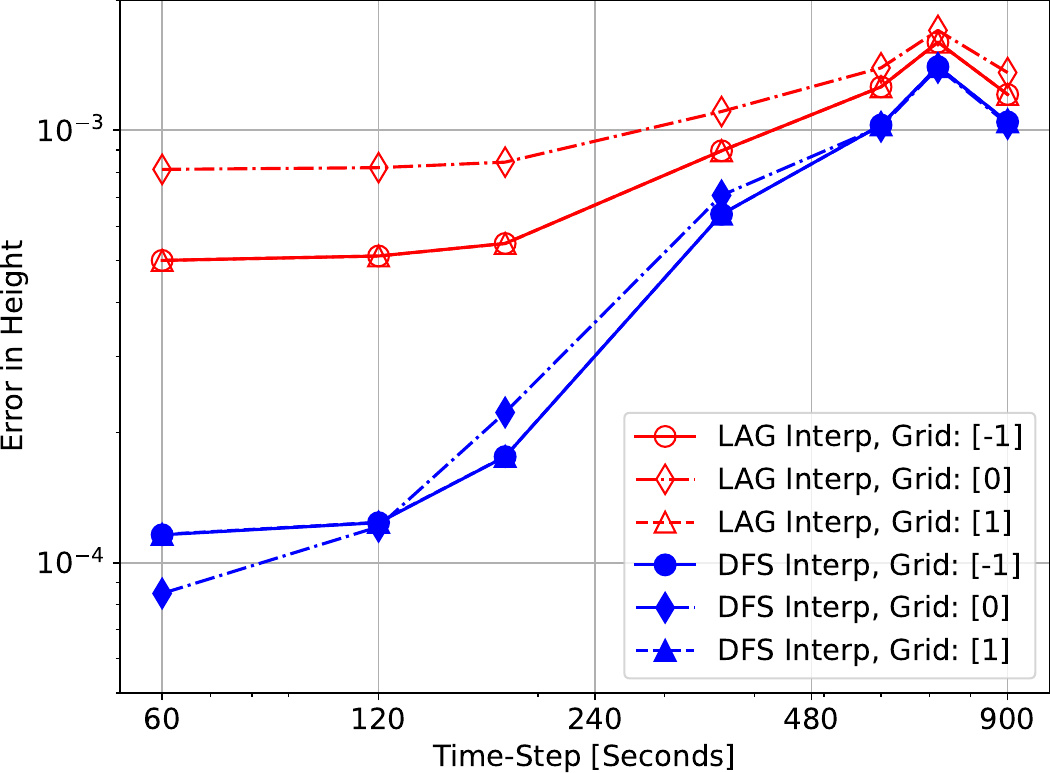}
\caption{LAG interpolation reference}
\label{fig:sw6tlag}
\end{subfigure}
\hfill
\begin{subfigure}[b]{0.45\textwidth}
\centering
\includegraphics[width=\textwidth]{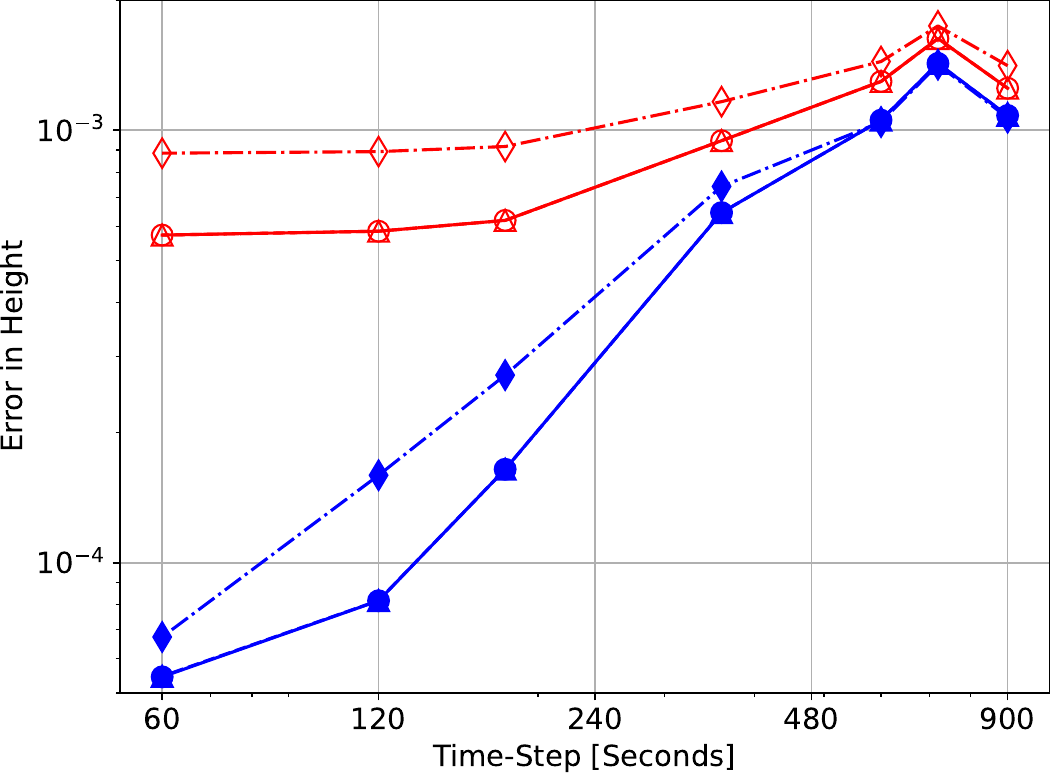}
\caption{DFS interpolation reference}
\label{fig:sw6tdfs}
\end{subfigure} 
\caption{Williamson TC6. Relative $L_2$ errors for the height field after 14 days of integration as the time step is decreased.  All results are for a grid size with $J=320$. Errors are computed using high-resolution SHT SISL solutions with (a) LAG and (b) DFS interpolation.\label{fig:sw6err}}
\end{figure}

Results on the temporal convergence of the height field are given in Figures \ref{fig:sw6tlag} and \ref{fig:sw6tdfs} using the LAG and DFS reference solutions, respectively. For the LAG case, the truncation wavenumber was set to $N=J-2$ for all time steps, but for the DFS case, it was set to $N=J-2$ for time step sizes less than or equal to 360 seconds and $N=J-8$ for larger time step sizes. As noted in~\cite{galewsky2004initial}, different truncation levels are often needed to avoid inducing an instability in this test case. The results in Figure~\ref{fig:sw6err} show that errors in height fields using the DFS interpolation converge at a much faster rate with decreasing $\Delta t$ compared to LAG interpolation, where convergence appears to stagnate for $\Delta t \leq 200$ s.  Additionally, the DFS interpolation shows better convergence properties when using the high-resolution DFS interpolation reference solution, whereas for LAG interpolation, the results are similar between the two reference solutions.  

\begin{figure}[htb]
\centering
\begin{subfigure}[b]{0.45\textwidth}
\centering
\includegraphics[width=\textwidth]{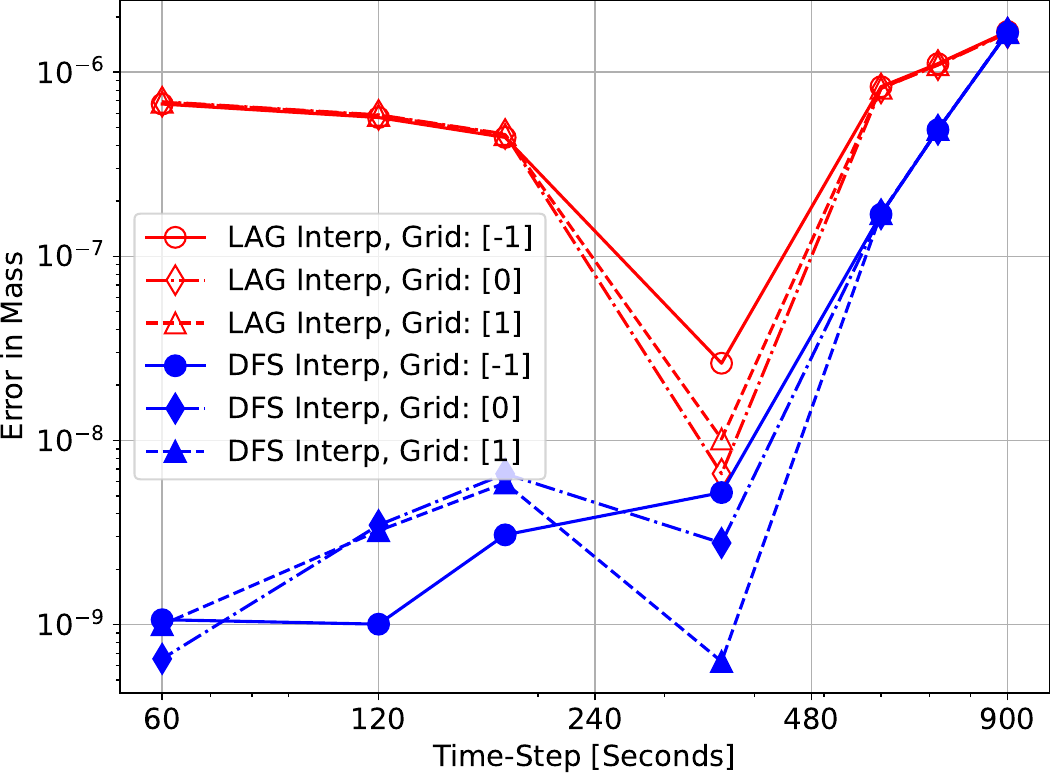}
\caption{Mass}
\label{fig:sw6masst}
\end{subfigure}
\hfill
\begin{subfigure}[b]{0.45\textwidth}
\centering
\includegraphics[width=\textwidth]{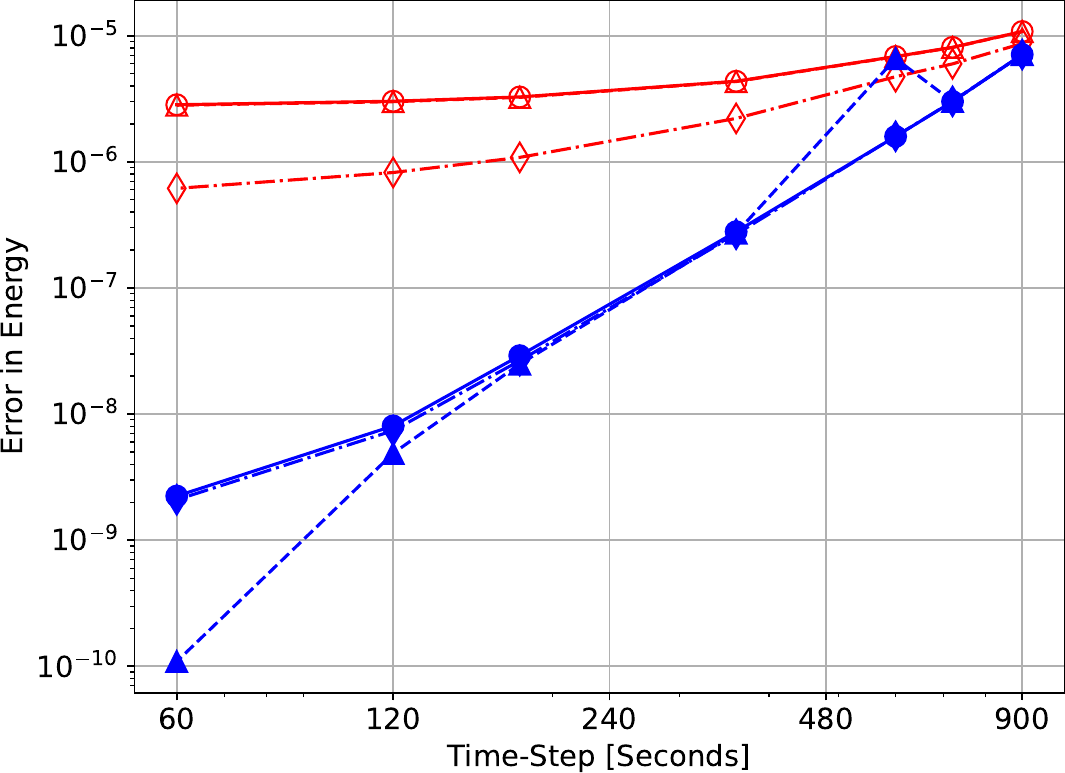}
\caption{Energy}
\label{fig:sw6energyt}
\end{subfigure}
\caption{Williamson TC6. Similar to Figure \ref{fig:sw6err}, but for errors in (a) mass and (b) energy after 14 days integration of TC6.}
\label{fig:tc6mass}
\end{figure}

Figures~\ref{fig:sw6masst} and~\ref{fig:sw6energyt} show the errors in the conservation of mass and energy of the methods as the time step decreases. Similar to TC5, we see that the DFS method does much better at conserving these quantities than the LAG method.

\begin{figure}[htb]
\centering
\begin{subfigure}[b]{0.45\textwidth}
\centering
\includegraphics[width=\textwidth]{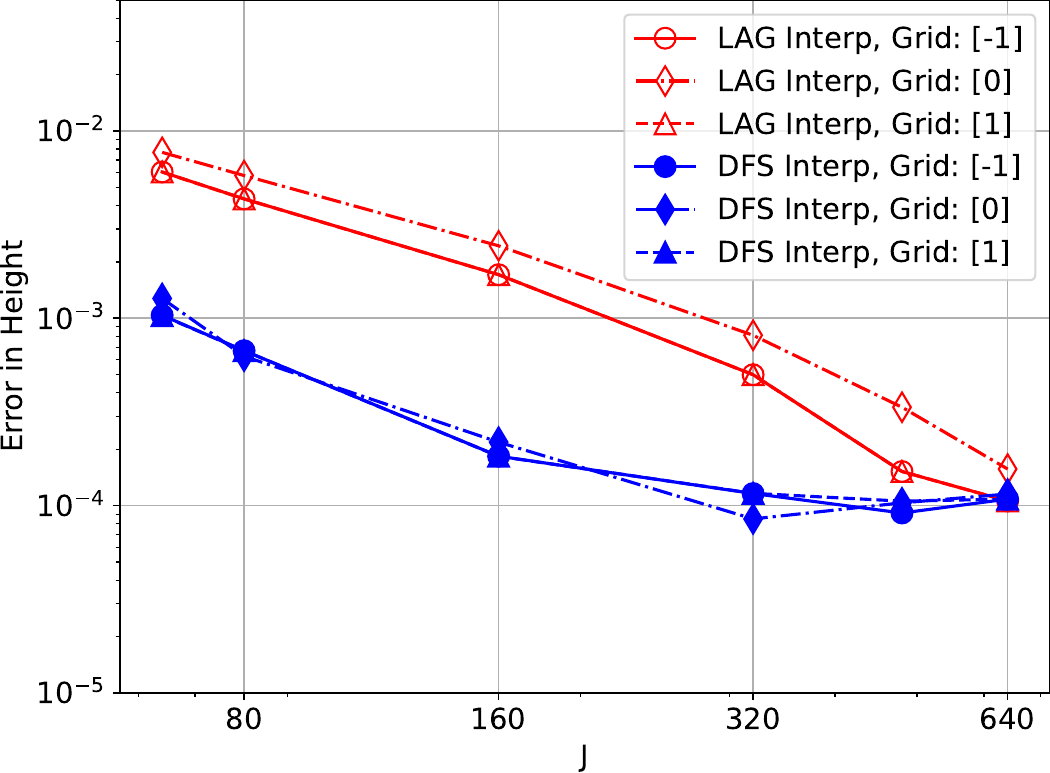}
\caption{LAG interpolation reference}
\label{fig:sw6splag}
\end{subfigure}
\hfill
\begin{subfigure}[b]{0.45\textwidth}
\centering
\includegraphics[width=\textwidth]{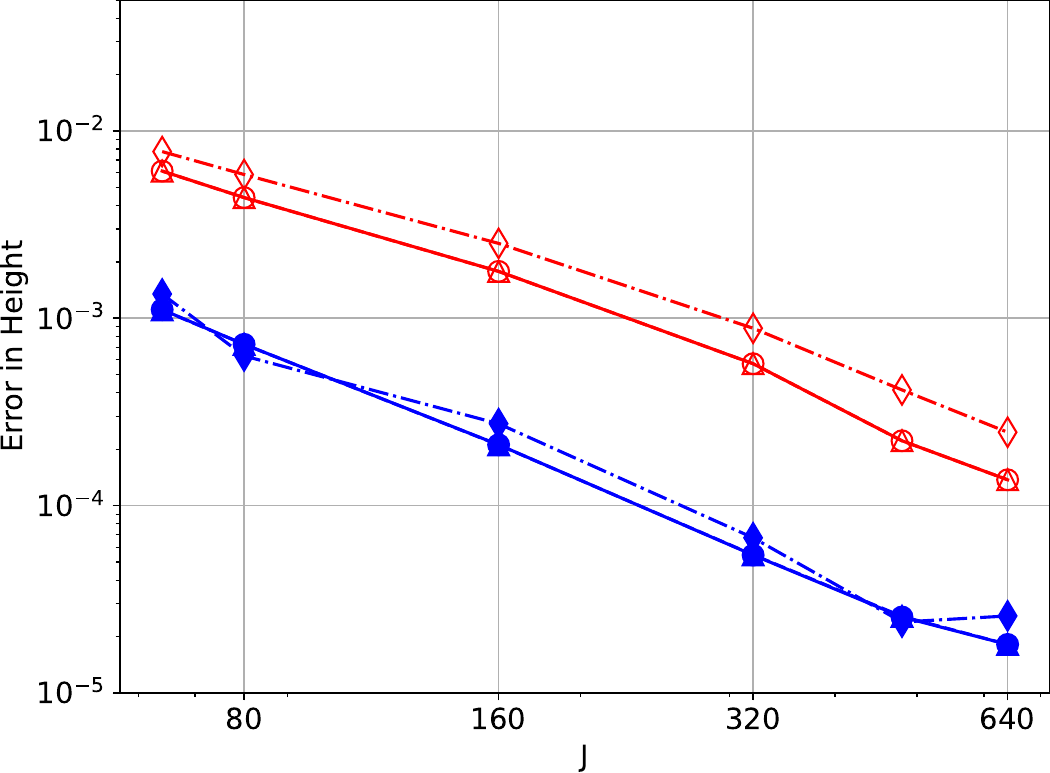}
\caption{DFS interpolation reference}
\label{fig:sw6spdfs}
\end{subfigure} 
\caption{Williamson TC6. Relative $L_2$ errors for the height field after 14 days of integration with a time step of 60 s and increasing spatial resolutions. Errors computed using high-resolution SHT SISL solutions with (a) LAG  and (b) DFS interpolation.\label{fig:sw6sp}}
\end{figure}

The spatial convergence of the two interpolation methods is shown in Figure \ref{fig:sw6sp} using the different reference solutions. For the LAG method, the truncation wavenumber is set to $N=J-2$, and for the DFS method, the truncation wavenumber is $J-2$ for $J<320$, $J-4$ for $J=320$, and $J-8$ for $J=480$ and $J=640$.  For both methods, we see that the $L_2$ norms of the errors in height fields appear to converge at the same rate with increasing $J$ when using the high resolution DFS reference.  However, the errors for the DFS method are about an order of magnitude lower for $J\geq 320$. For the LAG reference solution, convergence for the DFS method stalls for $J\geq 320$, indicating that the accuracy of the reference may be around $10^{-4}$.  This level of accuracy has also been noted in other spectral SWE models~\cite{spotz1998fast}.

\subsection{Galewsky Test Case}
\begin{figure}[htb]
\centering
\begin{subfigure}[b]{0.48\textwidth}
\centering
\includegraphics[width=\textwidth]{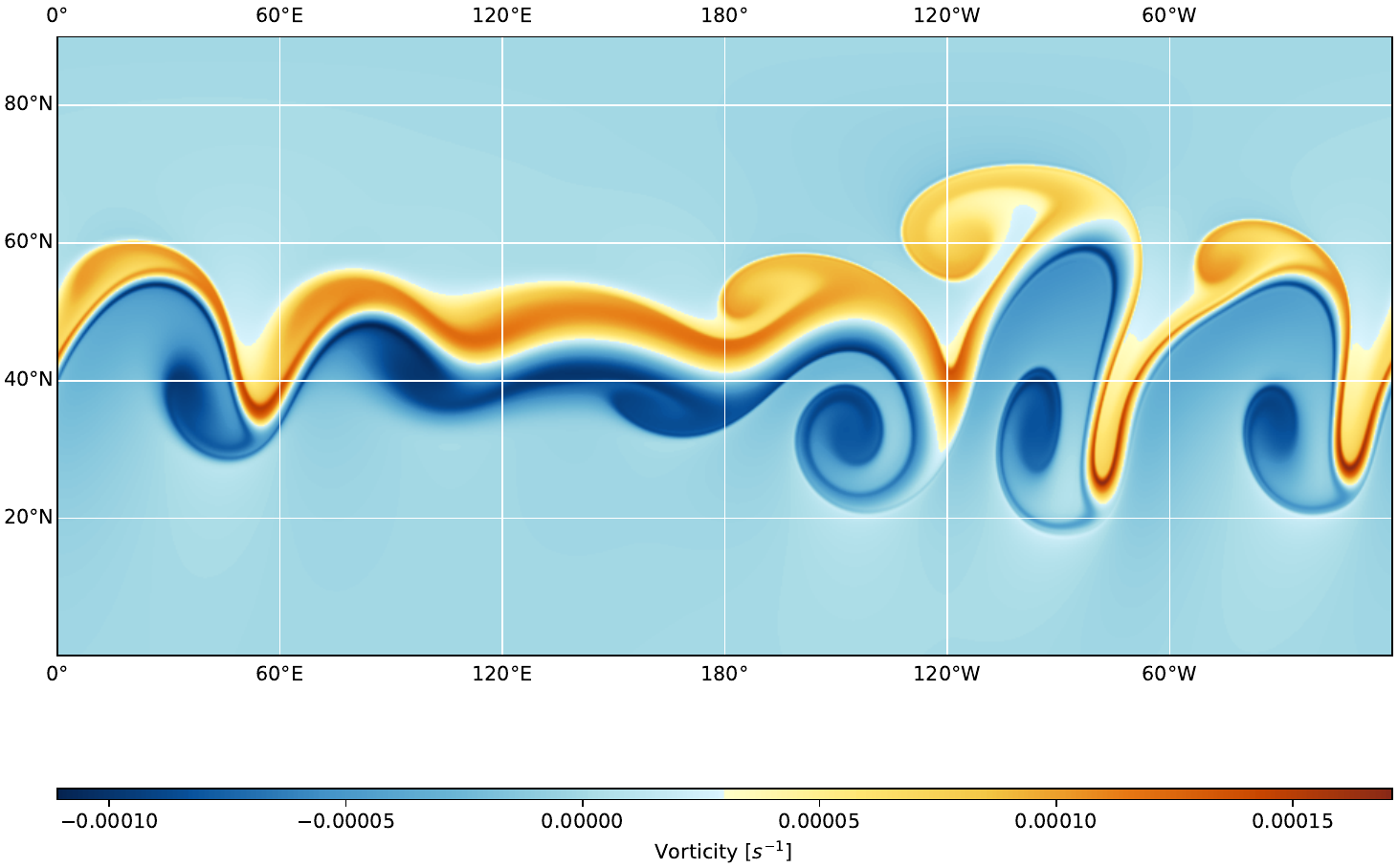}
\caption{LAG Interp}
\label{fig:sw10vortys}
\end{subfigure}
\hfill
\begin{subfigure}[b]{0.48\textwidth}
\centering
\includegraphics[width=\textwidth]{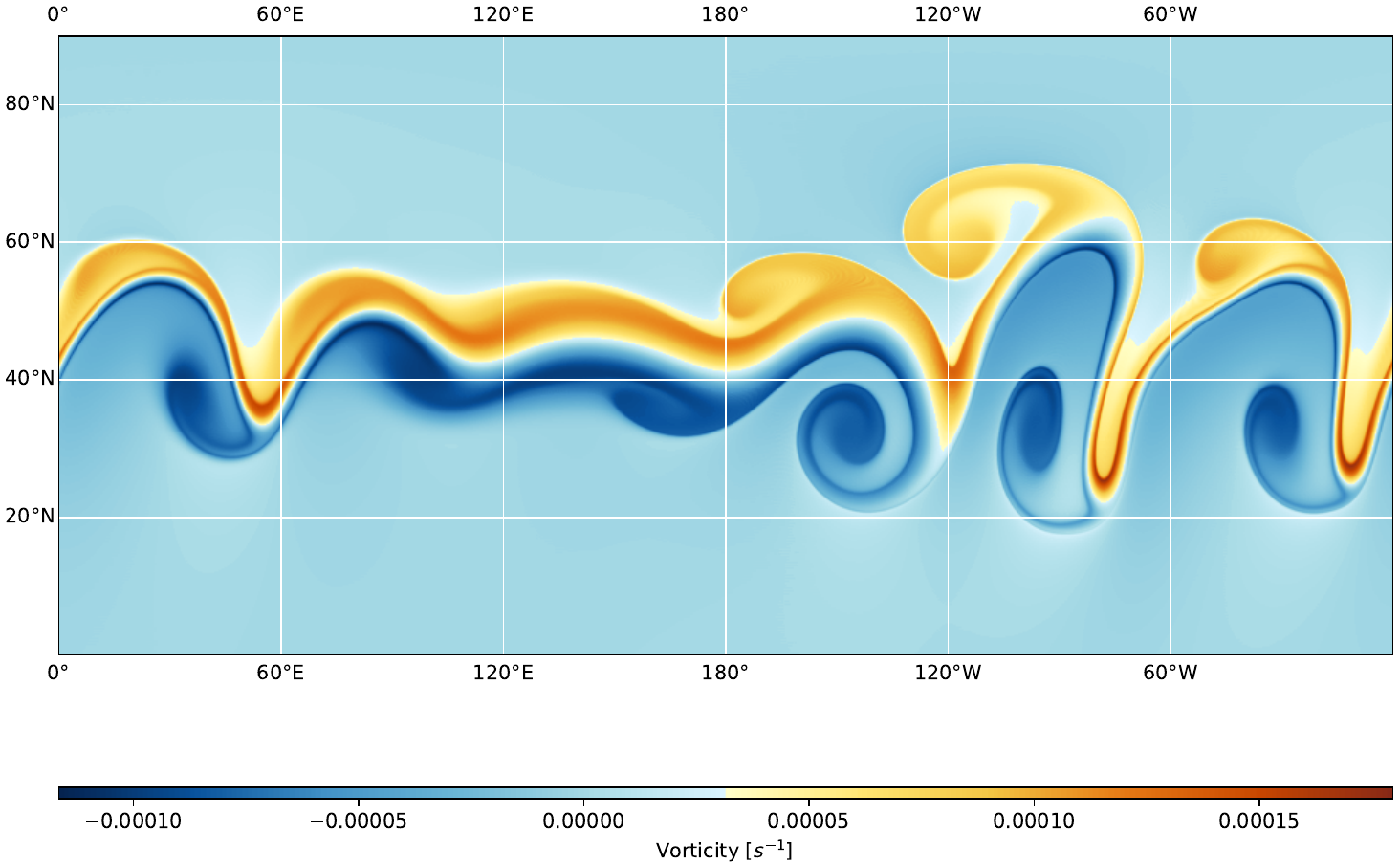}
\caption{DFS Interp}
\label{fig:sw10vort}
\end{subfigure} 
\caption{Galewsky test case. Predicted vorticity after 6-day integration for the two interpolation methods. Solutions for both methods were computed using grid [0] with $J=960$, a truncation wavenumber $N=J-2$, and time step $300$ s.\label{fig:sw10v}}
\end{figure}

The Galewsky test case models the evolution of a strongly nonlinear shallow-water flow, characterized by a rapid transfer of energy from large to small scales within a short timeframe~\cite{galewsky2004initial}. A localized perturbation to a zonal jet first generates high-frequency gravity waves that propagate around the sphere. These waves then roll up the flow into large, coherent vortices with sharp potential-vorticity gradients. By day 6, the vortices are fully developed and remain coherent, but later  their filaments begin to stretch and distort, and the flow gradually transitions toward more disordered, turbulent behavior. Figure \ref{fig:sw10v} shows the vorticity at day 6 computed with both the LAG and DFS methods, a stage dominated by these coherent vortex structures. At this time the solutions produced by the two methods are visually indistinguishable, and neither exhibits Gibbs-type oscillations near the sharp gradients.


\begin{figure}[htb]
\centering
\begin{tabular}{ccc}
\includegraphics[width=0.31\textwidth]{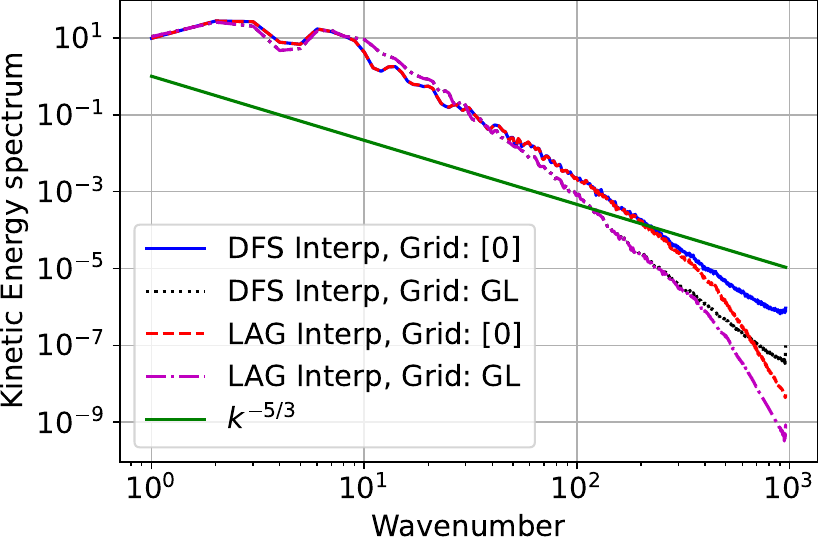} &
\includegraphics[width=0.31\textwidth]{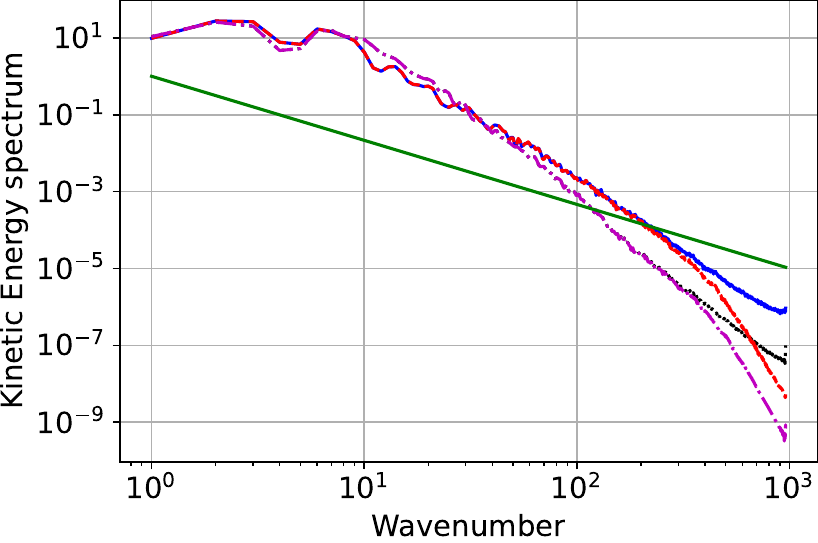} &
\includegraphics[width=0.31\textwidth]{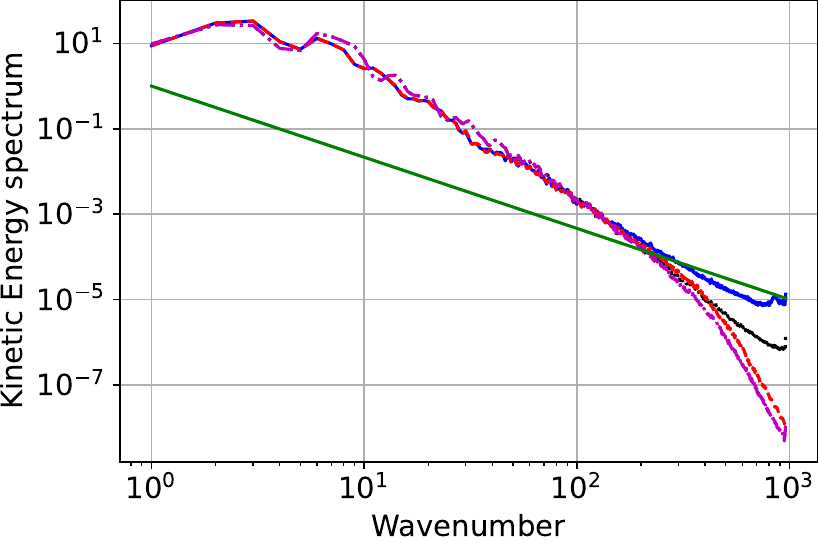} \\
(a) Day 6 & (b) Day 9 & (c) Day 12
\end{tabular}
\caption{
Galewsky test case. Kinetic energy spectrum of the horizontal winds $(\mathrm{m}^2/\mathrm{s}^2)$ at 6, 9, and 12 days of the simulation.  All solutions  were computed using $J=960$, a truncation wavenumber $N=J-2$, and time step $300$ s. \label{fig:sw10spectra}}
\end{figure}

Direct comparisons against a reference solution, as used for TC5 and TC6, become less informative for the Galewsky test because the flow evolves into a strongly nonlinear and eventually turbulent state for which no high-accuracy reference solution is available. Instead, it is standard to compare the statistical properties of the flow, such as the kinetic energy spectrum, to assess how well different numerical methods capture the distribution of energy across scales~\cite{wang2019understanding}. Figure~\ref{fig:sw10spectra} shows the kinetic energy spectrum of the horizontal winds $(\mathrm{m}^2/\mathrm{s}^2)$, as defined in \cite{lambert1984global}, after 6, 9, and 12 days of integration using four schemes: DFS SISL with LAG and DFS interpolation, and SHT SISL with LAG and DFS interpolation (denoted “Grid GL” in the figure). At day 6, when the vortices remain intact, the expected cascade of energy from low to higher wavenumbers is evident. By days 9 and 12, after the vortex structures begin to break down, the energy levels at low wavenumbers remain nearly identical across all four methods, indicating that each scheme captures the large-scale dynamics effectively. Differences among the schemes emerge primarily at the higher wavenumbers.

As time progresses and the flow becomes more turbulent, observational and modeling studies of rotating, stratified atmospheric flows suggest that the high-wavenumber horizontal kinetic energy spectrum tends toward a $k^{-5/3}$ scaling, where $k$ is the wavenumber~\cite{koshyk2001horizontal}. We observe in panels (b) \& (c) of the figure that the energy from both DFS-interpolation-based schemes does indeed approach the $k^{-5/3}$ behavior as time advances.  In contrast, the energy at high wavenumbers for both LAG-interpolation schemes decreases much more rapidly, indicating substantially stronger dissipation of small-scale energy.


We note that all methods exhibit oscillations in the kinetic energy spectrum at high wavenumbers, consistent with the observations reported in~\cite{yoshimura2022improved}. That study suggested that these oscillations could be reduced by using vorticity and horizontal divergence, rather than velocity, as prognostic variables. Evaluating this hypothesis lies beyond the scope of the present work and has not been examined here.

\subsection{Computational performance\label{sec:ctime}}
We conclude the numerical results with a comparison of the efficiency of the two interpolation methods, measured by the wall-clock time required to reach a comparable level of accuracy. For this evaluation, we record the total wall-clock time needed to achieve a prescribed error tolerance in the height field. This timing includes all components of the computation: semi-Lagrangian trajectory calculations, spectral transforms, and interpolation. The first two components use the same code configuration for both the LAG and DFS methods, so any differences arise entirely from the interpolation stage. All experiments were performed on a Linux workstation equipped with an Intel i9-9900X 3.5 GHz processor. OpenMP parallelization was used for the interpolation step in the LAG method and for the NUFFT computations in the DFS method, as supported by the FINUFFT library~\cite{Barnett2025}. The NUFFT allows the user to specify a tolerance parameter, where smaller tolerances produce more accurate results at a higher computational cost.  We report results for tolerances of $10^{-14}$ and $10^{-12}$, since larger tolerances did not provide sufficient accuracy.

\begin{figure}[htb]
    \centering
    \begin{subfigure}[b]{0.32\textwidth}
        \includegraphics[width=\textwidth]{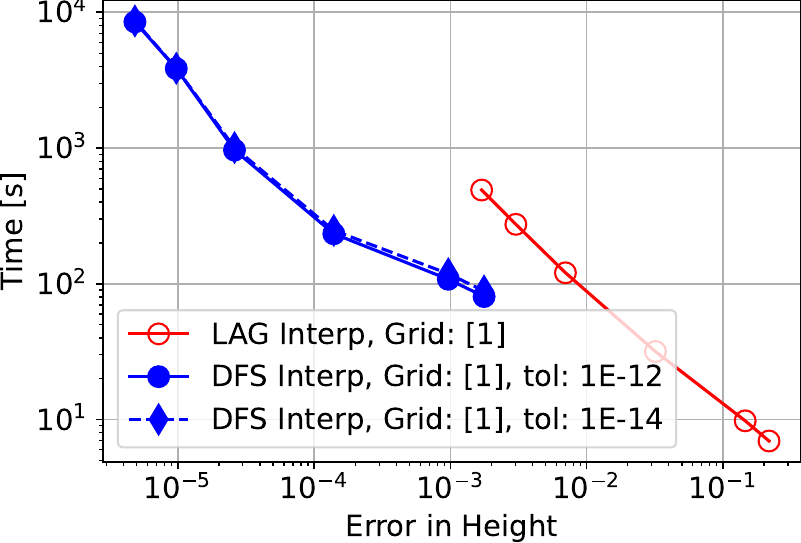}
        \caption{TC 1}
        \label{fig:TC1time}
    \end{subfigure}
    \hfill
    \begin{subfigure}[b]{0.32\textwidth}
        \includegraphics[width=\textwidth]{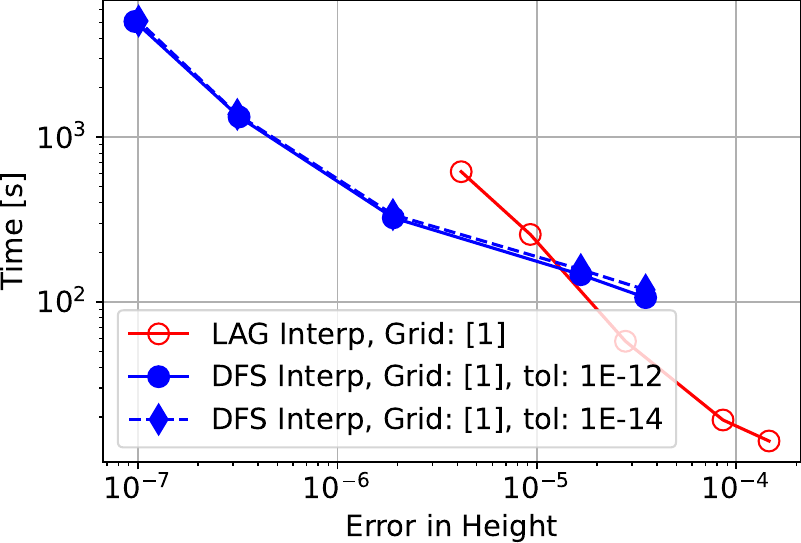}
        \caption{TC 5}
        \label{fig:TC5time}
    \end{subfigure}
    \hfill 
    \begin{subfigure}[b]{0.32\textwidth}
        \includegraphics[width=\textwidth]{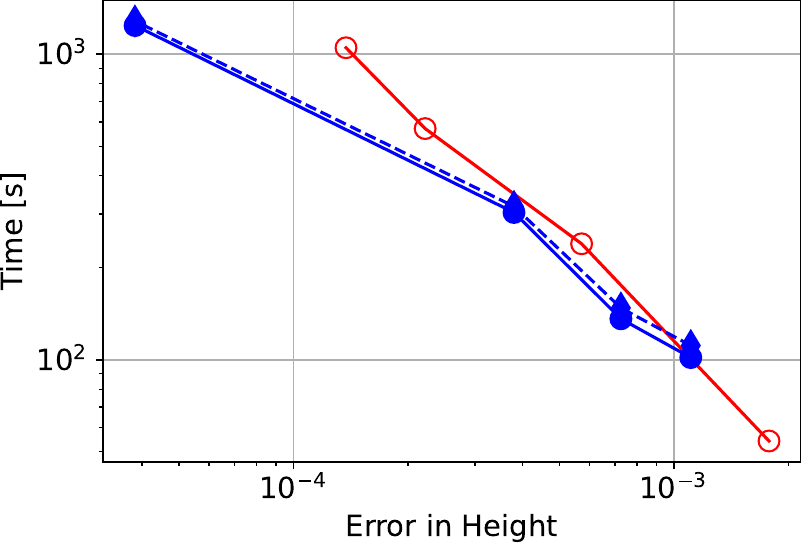}
        \caption{TC 6}
        \label{fig:TC6time}
    \end{subfigure}
    \caption{Comparison of the computational efficiency (wall-clock time versus error in the height field) for TC1, TC5, and TC6 using LAG and DFS interpolation. The DFS method includes results for NUFFT tolerances of $10^{-12}$ and $10^{-14}$.}
    \label{fig:comptime}
\end{figure}

Figure~\ref{fig:comptime}(a)–(c) presents the timing results for TC1, TC5, and TC6. For TC5 and TC6, the errors are measured against the high-resolution reference solutions described in Sections~\ref{sec:tc5} and \ref{sec:tc6}. The results show a consistent trend. For lower accuracy thresholds, the LAG interpolation method reaches the target tolerance more quickly. When higher accuracy is required, the DFS interpolation method becomes more efficient and ultimately reaches the target accuracy in less time. This demonstrates that DFS interpolation provides a computational advantage in regimes where high-order accuracy is important. It is also worth noting that accuracy is only one measure of performance. Results in the previous section show that DFS interpolation provides better conservation of key invariants and yields more physically realistic energy-cascade behavior than LAG interpolation. These findings offer additional evidence of the potential advantages of the DFS approach.


\section{Concluding remarks\label{sec:conclusion}}

Semi-implicit semi-Lagrangian (SISL) schemes enable stable integration of the shallow water equations (SWE) using time steps that are significantly larger than those permitted by corresponding Eulerian schemes, while still achieving comparable levels of accuracy. A common limitation, however, is that SISL schemes typically rely on relatively low-order interpolation during the advection step, a choice often viewed as a practical compromise between accuracy and computational efficiency. Although efficient, such interpolations can introduce excessive numerical damping, leading to a gradual loss of accuracy in long-term or highly nonlinear simulations.

In this study, we addressed this challenge by introducing a new spectrally accurate interpolation method based on the DFS representation and accelerated by the NUFFT to obtain near-optimal computational complexity. We incorporated this interpolation approach into both the DFS SISL scheme of~\cite{yoshimura2022improved} and the SHT SISL scheme of~\cite{yukimoto2019meteorological}. A comprehensive set of numerical experiments comparing the new spectral interpolation to standard Lagrange interpolation shows that the proposed method yields more accurate prognostic fields, demonstrates improved temporal convergence, enhances mass and energy conservation, and more faithfully captures the small-scale energy cascade in fully developed flows. In addition to these accuracy and conservation benefits, the DFS interpolation method also exhibits better computational performance for high-accuracy solutions, reaching stringent error tolerances in less wall-clock time than the Lagrange approach. These results motivate future work on extending DFS-based interpolation to more complex atmospheric models and fully three-dimensional dynamical cores.

\section*{Acknowledgements}
The authors thank Alex Townsend for early discussions on incorporating DFS interpolation into semi-Lagrangian schemes for the sphere and using the NUFFT to make it efficient. 

\paragraph{Funding}
The work of MC and GBW was partially supported by US National Science Foundation grants 2309712 and 2505987. The Flatiron Institute is a division of the Simons Foundation.

\bibliographystyle{siam}
\bibliography{swe}

\section*{Statements and Declarations}

\paragraph{Competing Interests}

The authors declare that they have no known competing financial interests or personal relationships that could have appeared to influence the work reported in this paper.

\paragraph{Author Contributions}

\begin{itemize}
\item MC: Conceptualization, Methodology, Software, Formal analysis, Investigation, Writing: Original draft preparation, Writing: Review \& Editing.

\item DF: Conceptualization, Methodology, Software, Resources, Writing: Review \& Editing, Supervision

\item GBW: Conceptualization, Methodology, Writing: Original draft preparation, Writing: Review \& Editing, Funding acquisition, Supervision, Software

\end{itemize}

\paragraph{Data Availability}

The software for reproducing the results in this study can be accessed at \url{https://github.com/michaelchiwere1/dfs-sisl-swe}.

\end{document}